%% file: MedinaPeszynska_2017_revised.tex
\documentclass{article}

\usepackage{lineno}
\usepackage{amsmath, amssymb,url}
\usepackage{graphicx,color}
\modulolinenumbers[5]

\input{amsnumbers.tex}

\newcommand\myskip[1]{}
\newcommand\mpcite[2]{\cite{#1}, (#2)}

\newcommand\norm[1]{{\parallel\! #1\! \parallel}}
\newcommand\abs[1]{{| #1|}}

\newtheorem{cor}{Corollary}
\newtheorem{lemma}{Lemma}
\newtheorem{proposition}{Proposition}
\newtheorem{remark}{Remark}
\def\ba{\begin{eqnarray}}
\def\ea{\end{eqnarray}}
\def\bas{\begin{eqnarray*}}
\def\eas{\end{eqnarray*}}
\def\R{\mathbf{R}}
\def\Fmat{\mathcal{F}}
\def\Bmat{\mathcal{B}}
\def\Cmat{\mathcal{C}}
\def\Dmat{\mathcal{D}}
\def\Lmat{\mathcal{L}}
\def\Amat{\mathcal{A}}
\def\Mmat{\mathcal{M}}
\def\tw{\tilde{w}}
\def\tW{\tilde{W}}

\def\tU{\tilde{U}}

\def\tB{\tilde{\Bmat}}
\def\tC{\tilde{\Cmat}}
\def\tF{\tilde{F}}

\newcommand{\mytitlenote}{Research partially supported by NSF
  DMS-1115827 ``Hybrid modeling in porous media'', and NSF DMS-1522734
  ``Phase transitions in porous media across multiple scales''}

\title{Stability for implicit-explicit schemes for non-equilibrium
  kinetic systems\\ in weighted spaces with
  symmetrization\footnote{\mytitlenote}}
	
\newcommand{\mycorrespondingauthor}{
Corresponding author\newline
URL: {https://www.wpi.edu/people/faculty/fpmedina} (F.~Patricia Medina)\newline
\hspace*{9mm}{http://www.math.oregonstate.edu/~mpesz} (Malgorzata Peszynska)
}

\author{F.~Patricia Medina \\
{\footnotesize \it Mathematical Sciences Department, 
Worcester Polytechnic Institute,
100 Institute Road,}\\
{\footnotesize \it Worcester, MA 01609-2980} \\[2mm]
Malgorzata Peszynska\footnote{\mycorrespondingauthor} \\
{\footnotesize \it Department of Mathematics, Oregon State University, Corvallis, OR 97331}
}

\date{}
\begin{document}
\maketitle

\begin{abstract}
We consider kinetic systems, and prove their stability working in
weighted spaces in which the systems are symmetric. We prove stability
for various explicit and implicit semi-discrete and fully discrete
schemes. The applications include advective and diffusive transport
coupled to the accumulation of immobile components governed by
non-equilibrium relationships. We also discuss extensions to nonlinear
relationships and multiple species.
\end{abstract}
	
	

\section{Introduction}
\label{sec:intro}

In this paper we transform and analyze semi-implicit numerical schemes
for an evolution system
\begin{subequations}
\label{eq:sys}
\ba
(\phi u)_t + v_t + \nabla \cdot(q u) - \nabla \cdot(\phi d \nabla u)=f,
\\
v_t=\alpha(g(u)-v) 
\ea
\end{subequations}
which arises in a variety of important applications, e.g., transport
in porous media with adsorption. Here $\alpha>0$ and $g(\cdot)$ is
monotone, with details below. The positive coefficient $\phi$
is the porosity.

For this system there is no maximum principle, and if $f=0$, there is
not even a natural conservation or stability principle in the natural
norms of $(u,v)$. Further, the analysis of the simple finite
discretization schemes with well known truncation errors, even when
$g$ is linear, has to deal with nonnormality, and is unnecessarily
complex, even when $g(\cdot)$ is linear.

The transformation we propose involves symmetrization, rescaling, and
a change of variables. Equivalently, we work in weighted spaces. We
exploit the symmetrization to prove strong stability of the problem
and of the associated numerical schemes, from which the natural error
estimates follow. For fully implicit schemes the framework of
m-accretive operators reduces the stability analysis to the
verification that the operator is m-accretive. However, for
implicit-explicit schemes this is not sufficient, and we draw upon
Fourier analysis.

\paragraph{\bf Overview} 
For the linear case when $g(u)=cu$, with $c>0$, the abstract form of
\eqref{eq:sys} has the structure of a linear kinetic system
\begin{subequations}
\label{eq:ade0}
\begin{eqnarray}
	U'+V'+LU&=& F \label{eq:ade01}\\ V'+\alpha(V-cU)&=&
        0, \label{eq:ade02}
\end{eqnarray}
\end{subequations}
with the unknowns $U,V:(0,T] \to H \times H$, where $H$ is an
  appropriate Hilbert space to be defined, and the source term$F:[0,\infty) \in
  H$ is given. The linear transport operator $L$ is defined in the
  sequel, and we will require for $L$ to be m-accretive to get strong
  stability.

Our main technical objective is to study the stability of
\eqref{eq:ade0} and of one--step implicit and implicit-explicit
discrete schemes for \eqref{eq:ade0}
\begin{subequations}
\label{eq:dde}
\begin{eqnarray}
	\frac{U^{n}-U^{n-1}}{\tau} +
        \frac{V^{n}-V^{n-1}}{\tau}+LU^{n*}&=& F^n
\label{eq-dde1}
\\ 
\frac{V^{n}-V^{n-1}}{\tau} +\alpha(V^n-cU^n) &=&
        0. \label{eq-dde2}
\end{eqnarray}
\end{subequations}
which is solved at every time step $n=1,2,\ldots$ for the
approximations $U^n,V^n$ to $u(\cdot,t_n),v(\cdot,t_n)$. Here
  $F^n \approx F(t_n)$.  This one-step scheme is fully implicit if
$n*=n$. Other schemes arise for $n*\neq n$. The analysis of
\eqref{eq:dde} involves consideration of spatial discretization as
well as of time discretization. Our technique of symmetrization allows
to demonstrate strong stability of the schemes in a weighted space,
even though the original system \eqref{eq:ade0} has nonnormal
operators.

Extensions of \eqref{eq:ade0} to nonlinear systems and to systems with
multiple components will be also discussed.
\paragraph{\bf Motivation and context} The problem \eqref{eq:sys} comes from applications in 
subsurface modeling such as the transport of contaminant undergoing
adsorption, or coalbed methane reservoir simulation, but cover also a
variety of other applications. In those problems \eqref{eq:sys}
represents the conservation of mass of some chemical component, with
$u$ denoting the mobile concentration, and $v$ representing the
immobile component, while $g(\cdot)$ a general monotone (increasing)
function. We provide details on the applications in
Section~\ref{sec:applications}.

Numerical analysis of \eqref{eq:sys} with non-equilibrium kinetics was
given in \cite{BKneq97} for diffusion only, with focus on
non-Lipschitz $g(\cdot)$ important for liquid adsorption. In
\cite{BKneq95} Lagrangian techniques for advection with
non-equilibrium adsorption and in \cite{DvDW94} the Lagrangian
transport combined with Galerkin approximation to diffusion were
analyzed. In addition, in a sequence of papers devoted to the scalar
conservation laws with relaxation terms \cite{STW97} a problem similar
to \eqref{eq:model} but without diffusion is studied, and convergence
order of $O(\sqrt{h})$ is established. In turn, in \cite{P13NumPDE} we
studied the stability of schemes for a single equation analogue of
\eqref{eq:model} without diffusion and where $v$ was eliminated, and
in \cite{PSY15} we extended the analysis to cover the linear case with
diffusion. Furthermore, previous results on stability of schemes of
\eqref{eq:ade0} for the case of initial equilibrium were shown in
\cite{P95c,P13NumPDE,PSY15}.

Our approach in this paper provides a unified framework for the
analysis of a variety of explicit and implicit finite difference
schemes for the non-equilibrium advection-diffusion problems.  In
particular, it establishes strong stability as well as optimal error
estimates of order $O(h)$ or $O(h^2)$.

\paragraph{\bf Outline} In Section~\ref{sec:applications} we motivate the study of
\eqref{eq:ade0}, provide examples of $L$, and provide literature
review. In Section~\ref{sec:abstract} we describe the main idea of
symmetrization in the abstract setting leading to the stability of the
numerical schemes. In Section~\ref{sec:discrete} we provide concrete
examples of fully discrete schemes for \eqref{eq:ade0} and evaluate
their stability, and in Section~\ref{sec:simulations} we illustrate
the theory with numerical examples, and convergence studies. We close
in Section~\ref{sec:extensions}, where we outline extensions to the
nonlinear and multi-species case, and discuss future work.

\paragraph{\bf Notation and assumptions} Throughout the paper we assume that
$c>0, \alpha >0$; otherwise, the system is decoupled and trivial.  $I$
always denotes the identity operator or matrix, as is clear from the content. 

With the original variables in \eqref{eq:sys} denoted by $u(x,t)$ and
$v(x,t)$, we consider the vectors
$U(t)=u(\cdot,t)=\left(u(x,t)\right)_x \in H$. Each $U(t),V(t)$ lives
in a Hilbert space $H$, with the inner product denoted by $\langle
\cdot,\cdot\rangle_H$; we drop the subscript $H$ when it does not lead
to a confusion. The domain of an operator $L$ is denoted by $D(L)$,
and the time derivative $U'(t)$ or $\frac{dU}{dt}(t)$ generalizes the
partial derivative $\frac{\partial }{\partial t}$, and is defined in
an appropriate abstract setting, such as that developed in
\cite{Show-monotone}.

The vector $W=W(t)=[U(t),V(t)]^T$ lives in $H \times H$, which is
endowed either with the natural or weighted inner product, with
details below. We also consider new variables $\tW$ in appropriate
spaces. The (matrices of) operators on $W$ or $\tW$ in the product
space $H \times H$ are denoted similarly to those on $H$. In
particular, for $L:D(L) \subset H\to H$ we define
$\Lmat=\left[\begin{array}{cc}L&0\\0&0\end{array}\right]$ on $D(\Lmat)
= D(L) \times H \subset H \times H$. We define $\Amat$, $\Dmat$
analogously.

In discrete schemes, we consider uniform time stepping $t_n=n\tau$
with $n=0,1,\ldots$, and time step $\tau$. For fully discrete schemes,
we denote the spatial grid parameter by $h$ and consider finite
dimensional analogues of the operators such as $L_h$ for $L$. For the
unknown $u(x,t)$ we denote by $u^n(x) \approx u(x,t_n)$ its
semi-discrete in time approximation, by $U^n \in H$ the collection
$(u^n(x))_x \in H$. In turn, $U_h(t)=(u_j(t))_j\in H_h$ where $u_j(t)
\approx u(x_j,t)$ is the semi-discrete in space approximation, in a
discrete (usually finite dimensional) subspace $H_h\subsetneq H$ of
dimension dependent on $h$. Finally, fully discrete approximations
$U_h^n \approx U_h(t_n)$. Most of our results are formulated on $H$,
and are shown to apply on $H_h$.

\section{Motivation and literature review}
\label{sec:applications}

In this section we develop the applications which motivate
\eqref{eq:ade0} and provide details on its abstract setup, in
particular on the properties of $L$ as they follow for the special
cases of \eqref{eq:sys} under assumed boundary conditions. Our
presentation of the model follows the literature on coalbed methane
adsorption \cite{ShiMazumder08,Kovscek2007,JessenKovc08} where
\eqref{eq:sys} arises directly, and finite volumes are used; see also
our expository work in \cite{PIMA11}. In Section~\ref{sec:nonlocal} we
discuss a particular direction in which \eqref{eq:ade0} is reduced to
a single equation; this is discussed in \cite{P13NumPDE} under initial
equilibrium assumption.

\subsection{Applications}
\label{sec:model}
In porous media, the model of transport with adsorption \eqref{eq:sys}
describes the evolution of the concentration of a chemical. Consider
an open bounded region of flow $\Omega \in \R^k,k=1,2,3$, in which the
volumetric flux $q$, with $\nabla \cdot q=0$, is given; assume also the
porosity $\phi(x)>0$ and the (uniformly positive definite) diffusion
coefficient $d$ are known.  If the chemical is adsorbed in the porous
medium, the mass conservation must include also the rate of change of
the adsorbed immobile amount denoted by $v$. The mass conservation of
the chemical being transported by advection $A u = \nabla \cdot (qu)$
and diffusion--dispersoon $Du = - \nabla \cdot(d \nabla u)$ , with
adsorption term, is
\begin{subequations}
\ba 
\label{eq:model}
(\phi u)_t + v_t + Au + Du =f, \;\;\; x \in \Omega, \;\ t>0\ea
and it remains to specify the relation of $v$ to $u$. 

The equilibrium relationship $v = g(u)$ which can be used to complete
\eqref{eq:model} assumes that the time scale of transport is much
slower than that of the adsorption.
In turn, the non-equilibrium or kinetic model
\ba
\label{eq:noneq}
v_t + \alpha(v-g(u)) = 0
\ea
\end{subequations}
allows to treat the time scales of adsorption and of transport on par
with each other, with $\alpha>0$ denoting the rate of the process. As
$\alpha \to \infty$, it is expected that \eqref{eq:noneq} has
solutions close to the equilibrium. The linear relationship $g(u) =
cu$, with $c>0$ is what is assumed throughout most of this paper.

For well-posedness, we require appropriate boundary conditions on $u$
as well as initial conditions for both $u$ and $v$. 

\subsubsection{Abstract setting}
\label{sec:problem}
In the abstract form, the model \eqref{eq:model} and \eqref{eq:noneq},
upon absorbing nonessential constants in the definitions of $u$, are
written as \eqref{eq:ade0}, in which $L=D+A$ is the abstract
diffusion-advection transport operator, and where \eqref{eq:ade0} is
posed as a Cauchy problem in an appropriate function space. 

Consider $L:D(L)\to H$ in a Hilbert space $H$, with domain $D(L)$. We
recall that $L$ is accretive if $\langle LU,U \rangle \geq 0$ for any
$U \in D(L)$. Additionally, $L$ is m-accretive if $I+L$ is onto $H$
that is, for any $F \in H$ the problem $U+LU=F$ is solvable (from
accretiveness there follows the uniqueness of the solution). 

For an m-accretive $L$, the following results are well known; see,
e.g., \mpcite{Show-monotone}{Sec.I.4}. The dynamics of $U'(t)+LU=0,
U(0)\in H$ is governed by a linear contraction semigroup, so that, in
particular, $U(t) \in D(L)$. If $L$ is self-adjoint, additional
regularity and convergence properties follow. The
  nonhomomogeneous case of $U'(t)+LU=F$ requires that $F\in
  (C^1[0,\infty),H)$. See, e.g. \mpcite{Show-monotone}{Prop.4.1} and
    \mpcite{Showalter77}{Cor. 3.B}.

For the applications of \eqref{eq:sys} described in
Sec.~\ref{sec:model} we consider $H=L^2(\Omega)$ with the inner
product $\langle \psi,\xi \rangle =\int_{\Omega} \psi\xi$.  In the
case of periodic boundary conditions, without loss of generality, one
can use $\Omega=(0,1)^k$, but we only analyze $k=1$ case.  We will
recall the standard abstract results for $Au = - d\nabla^2u + \nabla
\cdot(q u)$, with weak rather than the classical (partial)
derivatives. The definition of $L$ and $D(L)$ accounts for the
boundary conditions. For details on this abstract setup see
\mpcite{Showalter77}{Ex. IV.2., p108} and
\mpcite{Show-monotone}{Prop. I. 4.2, p21}. For periodic case, see
\mpcite{Showalter77}{Ex IV.1, p107}), and for advection see
\mpcite{Showalter77}{Example IV.1}.

\begin{remark}
\label{rem:L}
(i) Let $L=D$, with $d>0$, and with homogeneous Dirichlet boundary
conditions imposed. We have $D(L) = H_0^1(\Omega) \cap H^2(\Omega)$, and $L$
is m-accretive self-adjoint.  
(ii) As in (i), but with homogenous
Neumann conditions, $D(L)=H^2(\Omega)$, $L$ is m-accretive, and
self-adjoint. 
(iii) As in (i), with periodic boundary conditions,
e.g.. when $\Omega=(0,1)$ we have $D(L)=\{\psi \in H^2(\Omega),
\psi(0)=\psi(1), \psi'(0)=\psi'(1)\}$. The operator $L$ is m-accretive
and self-adjoint.
(iv) Case $L=A$, with properly posed conditions on the inflow
boundary, or with periodic boundary conditions, e.g., for
$\Omega=(0,1)$, $D(L) = \{ \psi \in H^1(\Omega): \psi(0)=\psi(1)\}$;
The operator $D$ is m-accretive but not selfadjoint.
(v) Case $L=A+D$, and $d>0$, $q \neq 0$. With periodic b.c., $D(L)$
  is as in (iii), and $L$ is m-accretive but not selfadjoint.
\end{remark}

\subsection{Related models and previous work}

In models of transport with gas adsorption, \eqref{eq:noneq} allows to
account for subscale diffusion accompanying the overall transport; see
\cite{KingErtekin86,ShiDur03}. More general models in which $\alpha$
is a monotone operator can be used, e.g., to model hysteresis in
adsorption \cite{PS98} or non-equilibrium phase transitions
\cite{DiBSho82}. Further, non-equilibrium relation \eqref{eq:noneq} is
used to model transport in media with multiscale character, such as in
the classical Warren-Root and Barenblatt models of double porosity
\cite{WarRoo63,BZK60}; see also modeling and analysis in
\cite{ADH,Show93a,HornSho90}, and numerical analysis in
\cite{P95c,KP12}.

In previous work for nonlinear $g(\cdot)$ \cite{DvDW94} proved
a-priori error estimates for Lagrangian-Galerkin methods for
\eqref{eq:sys}, and in \cite{BKneq97} the analysis is for diffusion
only. Our paper handles the advection and diffusion problem together
for linear $g(\cdot)$, and handles the analysis as well as implicit
and explicit numerical schemes in the same framework. 

\subsection{Nonlocal formulation in $u$ under the initial equilibrium assumption}
\label{sec:nonlocal}

One can reformulate the coupled kinetic system \eqref{eq:ade01} as a
single equation with nonlocal in time terms. Recall Volterra
convolution integral term defined by $\beta \ast U = \int_0^t
\beta(t-s)U(s) ds$. We solve \eqref{eq:ade02} for $V(t)$ in terms of
$U(t)$ and substitute to \eqref{eq:ade01} to give the following
Volterra integro-differential equation
\begin{subequations}
\label{eq:decouple}
\ba 
\label{eq:nonlocal}
U'+\alpha U' \ast \beta + AU=F + \beta(t) (V(0)-c U(0)),
\ea
solved for $U$, where $\beta(t)=\alpha e^{-\alpha t}$. The variable
  $V$ can be recovered from $U$ by
\ba
\label{eq:c}
V(t)=e^{-\alpha t}V(0) + \int_0^t \alpha c U(s) e^{-\alpha (t-s)} ds.
\ea
\end{subequations}
Now we see that the second term on the right hand side of
\eqref{eq:nonlocal} acts like a source/sink term decreasing with $t$,
and it vanishes under the assumption of initial equilibrium
\ba
\label{eq:initeq}
V(0)-c U(0)=0.
\ea
The one-way coupling in \eqref{eq:decouple} focuses the attention on
$U$ while keeping track of the memory effects expressed by $U' \ast
\beta$. 

The effect of memory terms isolated from the source term can be
studied if \eqref{eq:initeq} is assumed.  This approach was followed
in \cite{P95c,P13NumPDE,PSY15}. Strong stability of $u$ for the
numerical schemes was proven for $L=D$ in \cite{P95c}, and for linear
or nonlinear advection operator $L=A$ in \cite{P13NumPDE}, where we
exploited the positivity of the kernel $\beta$. Even though we did not
prove it, the numerical results suggested that a maximum principle
holds for $U$.

If \eqref{eq:initeq} cannot be assumed, the positive source term in
\eqref{eq:nonlocal} can be expected to disturb the maximum principle
and/or stability.  Indeed, a simple example in Sec.~\ref{sec:0d}
readily demonstrates it.

As concerns numerical schemes, in \cite{P95c} the kernel was also
allowed to be weakly singular, e.g., $\beta = O(t^{-1/2})$, which
corresponds to subscale diffusion, i.e., the case when there is
diffusion in \eqref{eq:ade02}.  A more general case of nonlocal terms
and of operators $L$ was studied in \cite{PSY15} in which the
multiscale model derived in \cite{PS07}, but initial data was assumed
in equilibrium.  In \cite{P13NumPDE} we conjectured experimentally
that the presence of the memory term $u' \ast \beta$ would lead to an
increased regularity of the solution $u$ when $\beta$ was weakly
singular, but this effect appears weaker for the bounded $\beta$.

Beyond porous media, the effect of non-equilibrum (relaxation) such as
in \eqref{eq:noneq}, was studied, e.g., in \cite{TveWin}, and is an
important component of pseudo-parabolic models \cite{BohmShow85b}.

\section{Stability for the abstract symmetrized evolution system}
\label{sec:abstract}
We start by motivating the symmetrization and discussing the
properties and well-posedness of the symmetrized system in the
abstract form on some general Hilbert spaces $H$. 

First we re-arrange \eqref{eq:ade0} in an equivalent form 
\begin{subequations}
\label{eq:ade}
\begin{eqnarray}
	U'-\alpha(V-cU)+LU&=& F, \label{eq:ade1}
\\ 
V'+\alpha(V-cU)&=& 0. \label{eq:ade2}
\end{eqnarray}
\end{subequations}
This arrangement is similar to those used in multiscale models such as
the Warren-Root or Barenblatt models \cite{WarRoo63,BZK60}. For these
models however $c=1$; their analysis and numerics can be found, e.g.,
in \cite{KP12}. When $c\neq 1$, the analysis requires additional work.

In a vector-matrix form with $w=[u,v]^T$ we write \eqref{eq:ade}
\ba
\label{eq:bsys}
W'+\Bmat W= W' + \Cmat W + \Lmat W  = [F,0]^T,
\ea
with
\ba
\label{eq:matrices}
\Cmat=\alpha \left[ \begin{array}{cc} cI &-I\\-cI &
I\end{array}\right], \;\; \Lmat = \left[ \begin{array}{cc} L &0\\0& 0\end{array}\right].
\ea
This system of evolution equations is solved for $W(t)=[U(t),V(t)]^T
\in H \times H$. The space $H \times H$ is endowed with the natural
inner product $(\cdot,\cdot)_{H \times H}$ and norm $\norm{\cdot}_{H
  \times H}$ on the product space, where
$\langle[U,V]^T,[\phi,\psi]^T\rangle_W=\langle U,\phi \rangle_H+
\langle V,\psi\rangle _H$.

\paragraph{\bf Challenge} The system \eqref{eq:bsys} is linear, and thus is trivially well-posed,
e.g., if $H=\R^P$, $P\in \N$.  However, in Section~\ref{sec:0d} we
show with a simple example on $H=\R$, that the system \eqref{eq:bsys}
is not stable in $\norm{w}_{H \times H}$, even though the solutions to the
homogeneous problem eventually decay to $0$.

Unless $c=1$, the operator $\Cmat$ and $\Bmat$ are not self-adjoint
and nonnormal with respect to $\langle \cdot,\cdot\rangle _{H \times H}$. In
consequence, the analysis of the numerical schemes for \eqref{eq:bsys}
is quite complicated. Therefore, we consider a weighted inner product
on $H \times H$ or, equivalently, a change of variables. This idea
which we explain in Sections~\ref{sec:symmetry} and
\ref{sec:alternative}, makes the subsequent analysis of numerical schemes
fairly straightforward.

\subsection{Symmetrization and rescaling}
\label{sec:symmetry}
We propose to consider a weighted (scaled) inner product $\langle
\cdot, \cdot\rangle_H$ on $H \times H$, and an associated norm
$\norm{\cdot}_c$
\ba
\label{eq:qoi}
\langle [U,V]^T,[\phi,\psi]^T\rangle_c=c\langle U,\phi\rangle_H+\langle V,\psi\rangle_H,\;\;
\norm{[U,V]^T}_c^2 = c \norm{U}^2_{H} + \norm{V}_{H}^2. \ea
In other words, instead of $H \times H$ we consider the new (Hilbert)
space $W_c$ which is $H \times H$ endowed with
$\langle\cdot,\cdot\rangle_c$. In the new space $W_c$ we are able to
prove the stability of the evolution system, and of the appropriate
numerical schemes for the diffusion-advection examples.

The use of weighted inner product space can be interpreted as changing
variables from $U$ to $\tU=\sqrt{c}U$, since $\norm{[U,V]^T}_c^2=
\norm{ [\sqrt{c}U,V]^T }_{H \times H} ^2$. We also denote the change
of variables from $w$ to $\tw$ with
\ba
\label{eq:tw}
\tW=[\tU,V]^T=[\sqrt{c}U,V]^T.
\ea

To show how we exploit the space $W_c$, we rewrite \eqref{eq:ade} by
scaling the first component equation of \eqref{eq:ade} by $\sqrt{c}$,
and re-distributing the appropriate constants, by linearity of $L$. We
see that \eqref{eq:ade} is equivalent to 
\begin{subequations}
\label{eq:tsys}
\ba
\sqrt{c}U'-\alpha( \sqrt{c} V-c \sqrt{c} U)+L\sqrt{c}U&=& \sqrt{c}F \label{eq:cde1}
\\ 
V'+\alpha(V-\sqrt{c} \sqrt{c}U)&=&
        0. \label{eq:cde2}
\ea
\end{subequations}
where we have also used $cU=\sqrt{c} \sqrt{c}U$ in
\eqref{eq:cde2}. Rewriting in the new variables
\ba
\label{eq:tbsys}
\tW'+ \tB \tW = \tW' + \tC \tW  + \Lmat \tW= [\tF,0]^T
\ea
with $\tF=\sqrt{c}F$ and  the operators defined as  
\ba
\label{eq:tmatrices}
\tB = \tC + \Lmat,\;\; 
\tC=\alpha \left[\begin{array}{c c} c I & -\sqrt{c} I \\ 
-\sqrt{c}I & I \end{array}\right].	
\ea

\subsection{Well-posedness}
Now we complete the formal discussion of the well-posedness of
\eqref{eq:bsys}.  We see that $B:D(B) \to H \times H$ with $D(B)=D(L)
\times H$ is {dense in } $H \times H$. Similarly, $\tB:D(\tB)\to W_c$,
and simply $D(\tB)=D(B)=D(L) \times H$.

\begin{proposition}
\label{prop:wellposedness}
Let $L$ be m-accretive on $H$. Then the
operator $B$ is m-accretive on $W_c$. Equivalently, $\tB$ is m-accretive on $H \times H$. 
\end{proposition}
{\noindent \sc Proof: }(i) The proof follows from \eqref{eq:tmatrices} and \eqref{eq:tw} by the calculation
\begin{multline}
\label{eq:proposition}
\langle \tB\tW,\tW \rangle _{H \times H}
\\
= \langle [(L+c\alpha I)\tU-\alpha \sqrt{c}I V,-\alpha \sqrt{c}I\tU+\alpha I V]^T,[\tU,V]^T\rangle_{H \times H}
\\
=\langle L\tU+c\alpha\tU-\alpha\sqrt{c}V,\tU\rangle_H
+\langle -\alpha\sqrt{c}\tU+\alpha V,V\rangle_H
\\
=\langle L\tU,\tU \rangle_H +c\alpha \langle \tU,\tU \rangle_H
-2\alpha \sqrt{c}\langle V,\tU \rangle_H +\alpha\langle V,V\rangle_H =
\\
\langle L\tU,\tU\rangle_H+\alpha\norm{\sqrt{c}\tU-V}_{H \times H}^2
\geq \langle L\tU,\tU \rangle_H \geq 0
\end{multline}
where we have exploited the symmetry 
$\langle V,\tU \rangle_H=\langle \tU,V\rangle_H$ 
and completed the square. The last step followed since $L$ is accretive on $H$. 

Similarly, we have that
\begin{multline}
\langle \Bmat W,W \rangle_c=\langle \Bmat [{U},v]^T,[{U},V]^T\rangle_c
\\
= \langle [(L + c \alpha I)U -\alpha  V, 
-\alpha cU + \alpha V]^T,[U,V]^T \rangle_c
\\
= c \langle LU,U \rangle _H+c^2\alpha \langle U,U \rangle_H-2\alpha c \langle U,V \rangle_H
+\alpha \langle V,V \rangle_H
\\
=c \langle LU,U \rangle_H + \alpha \langle cU-V,cU-V \rangle_H\geq c \langle LU,U \rangle_H\geq 0.
\end{multline}

(ii) To show that $\Bmat$ is m-accretive, i.e., that $I+\Bmat$ is onto $W_c$ we
show how to solve the system $(I+\Bmat)W = \Fmat$ for any $\Fmat=[F,G] \in H
\times H$.  To this end, we consider the solution of the stationary
counterpart of \eqref{eq:ade}
\bas
U -\alpha(V-cU)+LU = F,
\\
V+\alpha(V-cU) = G.
\eas
(In our problem \eqref{eq:ade0} we have $G=0$ but it is easy
  to consider the general case.)  Solving the second equation for $V$
in terms of $U$, back-substituting to the first equation, and
$\alpha-\frac{\alpha^2}{1+\alpha}=\frac{\alpha}{1+\alpha}$, we see
that $U$ satisfies
\ba 
\left((1+\frac{c\alpha}{1+\alpha}c)I+ L\right)U = F + \frac{\alpha}{1+\alpha}G
\ea
which can be solved for any $F\in H$, since $L$ is m--accretive.

\begin{cor}
\label{cor:stability}
Assume $\alpha,c>0$ and $L:D(L)\to H$ is m-accretive, and $\tW_{init}
\in W_c$.  By Hille-Yosida Theorem as quoted in
\cite{Show-monotone}{Prop.4.2, p21} and \cite{Show-monotone}{Thm I.5.1,
  p25}, we conclude that there exists a unique solution to the Cauchy
problem, with $\tw(t) \in D(\tB)$
\bas \tW'+\tB
\tW=[\tF,0]^T,\;\; \tW(0)=\tW_{init} \in H. 
\eas
The evolution if $\tW$ is governed by the linear contraction
semigroup. When $\tF=0$, we have the stability
\ba
\label{eq:stability}
\frac{d}{dt} \norm{\tilde{W}}^2 \leq 0.
\ea
\end{cor}

\subsection{Alternative motivation for symmetrization}
\label{sec:alternative}
We provide here another way to motivate the symmetrization and
rescaling proposed in Section~\ref{sec:symmetry}. 
We consider the homogeneous case of \eqref{eq:ade}, and take the inner
product of each component equation with $u$ and $v$, respectively. We
obtain
\begin{eqnarray*}
\langle U',U\rangle-\alpha\langle V,U\rangle 
+\alpha \langle cU,U \rangle + \langle LU,U \rangle &=& 0 
\\
\langle V',V \rangle +\alpha \langle V,V \rangle -\alpha \langle cU,V \rangle &=&0.
\end{eqnarray*}
Adding these identities directly does not produce useful results for
stability in $\norm{(U,V)}$, because the cross-terms do not cancel.
However, up to the scaling, the second term in the first identity is
similar to the third one in the second identity.  Multiplying the
first equation with $c$, and adding the resulting equations, we obtain
\bas
c \langle U',U \rangle-\alpha \langle V,cU \rangle +\alpha \langle cU,cU \rangle 
+c \langle LU,U \rangle 
+
\langle V',V \rangle +\alpha \langle V,V\rangle -\alpha \langle cU,V \rangle &=&0
\eas
Rearranging the terms, by symmetry of the inner product, we get
\begin{multline}
\label{eq:step}
c \langle U',U \rangle  +\langle V',V\rangle  
+\alpha \langle V,V\rangle - 2 \alpha \langle V, cU \rangle 
\\+ \alpha  \langle cU,cU \rangle
+c \langle LU,U \rangle) = 0
\end{multline}
Next, for the first two terms in \eqref{eq:step} we write
\begin{eqnarray*}
c \langle U',U\rangle 
+ \langle V',V \rangle
= c \frac{1}{2} \frac{d}{dt} \norm{U}^2 + \frac{1}{2} \frac{d}{dt} \norm{V}^2 = 
\frac{1}{2} \frac{d}{dt} \norm{[\sqrt{c}U,V]^T}^2.
\end{eqnarray*}
The next three terms in \eqref{eq:step} are easily combined to give
$\alpha \langle V-cU,V-cU\rangle= \alpha \norm{V-cU}^2\geq 0$. Since $L$ is
accretive, upon $c \langle LU,U\rangle\geq 0$ we obtain from \eqref{eq:step}
\ba
\frac{d}{dt} \norm{[\sqrt{c}U,V]^T}^2 \leq 0.
\ea
In other words, we see that the system \eqref{eq:ade} is stable in the
quantity of interest $\norm{[\sqrt{c}U,V]^T}$, or in $[U,V]_{W_c}$.

\section{Stability of numerical schemes}
\label{sec:discrete}

In this section we discuss numerical schemes for \eqref{eq:bsys} and
their stability and convergence properties. We focus on one-step time-discrete
schemes for \eqref{eq:bsys} solved for $W^n \approx
w(\cdot,t_n) \in H$
\ba
\label{eq:bdiscrete} 
\frac{1}{\tau} (W^n-W^{n-1})+\Cmat W^n + \Lmat W^{n*}=F^n, \;\; n \geq 1.
\ea
If $n*=n$, the scheme is fully implicit, and if $n*=n-1$, we have
implicit-explicit schemes. Note that our treatment of the (stiff)
coupling term $\Cmat W^n$ is always implicit. Here $F^n$ is some
appropriately defined time-discrete approximation to $F(\cdot,t_n)$,
and $W^0$ is known from the initial conditions. 

We also note, upon \eqref{eq:tsys}, \eqref{eq:tmatrices},
that \eqref{eq:bdiscrete} is equivalent to 
\ba
\label{eq:tbdiscrete} 
\frac{1}{\tau} (\tW^n-\tW^{n-1})+\tC \tW^n + \Lmat \tW^{n*}=\tF^n, \;\; n \geq 1,
\ea
where $\tC$ is symmetric.

In fully discrete schemes, the abstract operators $\Lmat,\Cmat$ are
replaced by their finite dimensional analogues $\Lmat_h,\Cmat_h$
depending on the spatial discretization parameter $h$, and they are
solved for the vectors of spatial unknowns $W_h^n=(w_j^n)_j$ where
$w_j^n \approx w(x_j,t_n)$ as in
\ba
\label{eq:bdiscrete-full} 
\frac{1}{\tau} (W_h^n-W_h^{n-1})+\Cmat_h W_h^n + \Lmat_h W_h^{n*}=F_h^n,
\ea
with an analogous version for \eqref{eq:tbdiscrete}, which we
skip. Here $F_h^n$ is an appropriate discretization of $F$.

In addition, we recall that finite element formulations lead, instead of
\eqref{eq:bdiscrete-full}, to
\ba
\label{eq:mbdiscrete-full} 
\frac{1}{\tau} \Mmat_h (W_h^n-W_h^{n-1})+\Cmat \Mmat_h W_h^n + \Lmat_h W_h^{n*}=F_h^n, \;\; n \geq 1.
\ea
where $\Mmat_h=\left[\begin{array}{cc}M_h&0\\0&M_h\end{array}\right]$,
and $M_h$ is the symmetric positive definite mass (Gram) matrix. For
generality, we adopt \eqref{eq:mbdiscrete-full} as the general fully
discrete formulation, since \eqref{eq:bdiscrete-full} is its special
case upon setting $M_h=I$.

\subsection{Fully implicit schemes for m-accretive $L$ and $L_h$}
We consider here $n=n^*$ in \eqref{eq:bdiscrete} or
\eqref{eq:mbdiscrete-full}. 

First observation is somewhat surprising. One might expect when
$F=0$, that
$ \norm{W^{n}}_{H \times H}\leq \norm{W^{n-1}}_{H \times H}$, but this
does not hold, e.g., if $\Bmat$ is nonnormal. (See
example in Sec.~\ref{sec:0d}).

However, based on the discussion in Sec.~\ref{sec:abstract}, stability
can be shown easily in weighted spaces.
\begin{lemma}
\label{lem:fi}
Let $L$ be m-accretive. Then the fully implicit scheme
\eqref{eq:tbdiscrete} is strongly stable in the weighted
spaces. We have, when $F=0$, that
\ba
\label{eq:sstability}
\norm{W^{n}}_c \leq \norm{W^{n-1}}_c \;\; \Longleftrightarrow
\norm{\tW^{n}} \leq \norm{\tW^{n-1}}, 
\ea
i.e., the operator $(I+\tau \tB)^{-1}$ is a contraction. For $F \neq 0$, we have 
\ba
\label{eq:fsstability}
 \norm{\tW^{n}} \leq \norm{\tW^{n-1}} + \norm{\tF^n}.
\ea
For \eqref{eq:mbdiscrete-full} we have 
\ba
\label{eq:msstability}
\norm{\Mmat_h^{1/2}W^{n}}_c \leq \norm{\Mmat_h^{1/2}W^{n-1}}_c \;\; \Longleftrightarrow 
\norm{\Mmat_h^{1/2}\tW^{n}} \leq \norm{\Mmat_h^{1/2}\tW^{n-1}} 
\ea
\end{lemma}

{\noindent \sc Proof: }The proof is immediate when we rewrite \eqref{eq:tbdiscrete}
as
\ba
\label{eq:tauB}
\frac{1}{\tau} (\tW^n-\tW^{n-1})+\tB W^n=\tF^n, \;\; n \geq 1.
\ea
Rearranging, and taking the inner product with $\tW^n$ we obtain 
\bas
\langle (I+\tau \tB) \tW^{n},\tW^n\rangle =\langle \tW^{n-1},\tW^n\rangle 
+ \tau \langle \tF^n,\tW^n\rangle . 
\eas
Since $L$ is accretive, by Proposition~\ref{prop:wellposedness} so is
$\tB$ on $H \times H$. Applying this property and Cauchy-Schwartz inequality we get
\begin{multline*}
\norm{\tW^n}^2 = \langle \tW^{n},\tW^n\rangle 
\leq 
\langle \tW^{n},\tW^n\rangle + \tau \langle \tB \tW^{n},\tW^n\rangle
\\
= \langle (I+\tau \tB) \tW^{n},\tW^n\rangle =\langle \tW^{n-1} +\tau F^n,\tW^n\rangle \leq 
(\norm{\tW^{n-1}}+\tau \norm{\tF^{n}})   \norm{\tW^n}.
\end{multline*}
where we have also used \eqref{eq:tauB}. 
Upon dividing by $\norm{\tW^{n-1}}$ we get
\eqref{eq:sstability}. Alternatively, we start from
\eqref{eq:bdiscrete} and work in weighted spaces in which $\Bmat$ is
accretive.

To prove \eqref{eq:msstability}, we proceed analogously using the
properties of $\tB_h$ on $H_h$ for
\eqref{eq:mbdiscrete-full}. Additionally we carry out an easy
calculation similar to \eqref{eq:proposition} involving $M_h$, which
takes advantage of positive definitiness and symmetry of $M_h$.

\begin{remark}
\label{rem:Lh}
Lemma~\ref{lem:fi} reduces the stability analysis of implicit schemes
to the verification whether $L$ (or $L_h$) is accretive. In addition,
for finite difference schemes the result \eqref{eq:fsstability}
applies directly to the error analysis, since \eqref{eq:bdiscrete} can
be interpreted as the error equation, in which the right-hand-side
represent the truncation error.  We see, in particular, that the error
accumulates linearly.    
\end{remark}

We collect several results for $L=D$ in Sec.~\ref{sec:D}; these are in
the framework of the method of lines (MOL). However, for some $L$ even
those covered in Rem.~\ref{rem:L}, and some schemes, $L_h$ is
non-symmetric, and it is hard to verify if it is m-accretive even if
$L$ is. In particular, for $L=A$ or $L=D+A$, and non-implicit schemes,
we proceed by von-Neumann analysis.

\subsection{FD discretization for diffusion with Dirichlet boundary conditions}
\label{sec:D}
First we consider FD discretization. For the sake of exposition, we
provide details for $k=1$ and $\Omega=(0,1)$.  We seek the interior
values $u_j^n$, $j=1,\ldots N_h$, where $h=\frac{1}{N_h+1}$ is the
spatial grid parameter. We also seek $v_j^n$ on the same grid of
interior points.

After symmetrization and rescaling, at every time step, one solves the
problem \eqref{eq:bdiscrete-full},
rewritten as
\begin{subequations}
\label{eq:sdde-dirichlet}
\begin{eqnarray}
	\frac{\tU_h^{n}-\tU_h^{n-1}}{\tau} -\alpha(\sqrt{c}V_h^n-c\tU^n)+L_h\tU_h^{n*}&=& \tF_h^n
\label{eq-sdde1-d}
\\ 
\frac{V_h^{n}-V_h^{n-1}}{\tau} +\alpha(V_h^n-\sqrt{c}\tU_h^n) &=&
        0. \label{eq-sdde2-d}
\end{eqnarray}
where the well known Dirichlet matrix is 
\ba
\label{eq:matrixD}
L_h:=\frac{d}{h^2} \left[\begin{array}{r r r r r r}  2 & -1 & 0 & \dots & 0 & 0
\\ -1 & 2 & -1 & 0 & \dots & 0 
\\ \vdots & \vdots & \ddots & \ddots & \ddots& 
\\ 0  & \dots & 0 & -1 & 2 & -1
\\ 0 & 0 & \dots & 0 & -1 & 2 \end{array} \right].
\ea
\end{subequations}
In addition, we note that $M_h$ in \eqref{eq:bdiscrete-full} is the identity matrix. 
Also, we recall that $L_h$ is symmetric positive definite on $\R^{N_h}$,
thus m-accretive. In fact, for this chosen domain $\Omega$, the eigenvalues of
$L_h$ are $\frac{d}{h^2} 2(1-cos(p\pi h))>0$, $p=1,\ldots N_h$.

When $n*=n$, \eqref{eq:sdde-dirichlet} has the form of
\eqref{eq:mbdiscrete-full}, with $M_h=I$. We can thus apply
Lemma~\ref{lem:fi}.  We obtain the following result which holds for
other domains $\Omega$, and Dirichlet boundary conditions, as long as
$L_h$ is symmetric and positive definite.

\begin{cor}
\label{cor:iD}
 The implicit in time scheme \eqref{eq:sdde-dirichlet} for $n*=n$ is
 strongly stable in the sense of \eqref{eq:msstability}.
\end{cor}

When $n*=n-1$, we consider the homogeneous version of
\eqref{eq:sdde-dirichlet} in the form
\ba
\label{eq:eD}
	\mathcal{H}_1 \tW^{n}=\mathcal{H}_0 \tW^{n-1},
\ea
where $\mathcal{H}_1$ and $\mathcal{H}_0$ are the block matrices, with
$b=\alpha \tau$,
\bas
\mathcal{H}_1=\left[
\begin{array}{c | c}
(1+bc)I & -b\sqrt{c}I\\ \hline
-b\sqrt{c}I	 & (1+b)I
\end{array}\right], \;\;\; 
\mathcal{H}_0=\left[
\begin{array}{c | c}
I-L_h & O\\ \hline
O	 & I
\end{array}\right].
\eas
Observe that the matrix $I-L_h$ is
symmetric, thus so is $\mathcal{H}_0$. If $X_p,p=1,\ldots P$ are the
eigenvctors for $L_h$, it is easy to show that each $[X_p,0]^T$ is an
eigenvector for $\mathcal{H}_0$, corresponding to the eigenvalues
$\lambda_p,p=1,\ldots P$ of $L_h$. In turn, the remaining eignevectors
of $\mathcal{H}_0$ are in the form of $[0,Y]^T$ where $Y\in\R^P$ is
arbitrary, with eigenvalue $\lambda=1$ of multiplicity $P$. The set of
eigenvalues for $\mathcal{H}_0$ is $\left\{1-\lambda_p\right\}.$ 
Since $\mathcal{H}_1$ is self-adjoint, conditional stability of
\eqref{eq:eD} follows by checking that the eigenvalues of
$\mathcal{H}_0$ are not exceeding $1$.

\begin{cor}
\label{cor:eD}
The scheme \eqref{eq:sdde-dirichlet} for $n*=n-1$ is
 conditionally strongly stable in the sense of \eqref{eq:msstability} if 
$\frac{2d\tau}{h^2} \leq 1 $. \end{cor}

Next, we briefly mention that handling periodic and Neumann problems
with FD is done differently than by constructing the simple analogue
of \eqref{eq:sdde-dirichlet}. Also, $L_h$ is typically not
symmetric. We do not discuss these cases here.

\subsection{FE discretization for $L=D$ and general boundary conditions}
Next we consider $L=D$ with Dirichlet, or Neumann, or periodic
boundary conditions covered by Rem.~\ref{rem:L} so that $L$ is
m-accretive on $D(L) \subset L^2(\Omega)$. Note that $L$ can have
variable coefficients and possibly correspond to some other boundary
conditions, as long as $L$ is m-accretive.

Next considered piecewise linear finite elements forming the
approximating subspace $V_h\subset H^1(\Omega)$, where $V_h$ accounts
properly for the essential boundary conditions. The nodal degrees of
freedom $u_j^n$, $v_j^n$, $j=1,\ldots N_h$ are in the space
$H_h=\R^{N_h}$. It is well known \cite{Thomee,Quarteroni} that the
matrix $L_h$ inherits the properties of the operator $L$ and, in
particular, it is symmetric and nonegative definite, thus
m-accretive. For this setting we have the result as follows. 

\begin{cor}
Finite element discretization \eqref{eq:mbdiscrete-full} for $n*=n$ is
strongly stable in the sense of \eqref{eq:msstability}.
\end{cor}

\subsection{Stability for advection and advection--diffusion via extended von-Neumann analysis for systems}
The von-Neumann framework for stability analysis of finite difference
schemes on uniform spatial grids for scalar linear equations with
constant coefficients is well known and is covered in various
textbooks see, e.g., \mpcite{LeVeque}{Chapters 9 and 10}. The
classical monograph \cite{MortonRicht} deals also with nonlinearity,
non-constant coefficients and coupled systems with non-normal
amplification matrix; we adopt their notation.

\paragraph{\bf Notation} 
First we establish the notation and recall the usual steps. We
consider the true solution $s(x,t) \in \R, x \in \R,t>0$, to a scalar
differential equation, which is approximated by a finite difference
equation with uniform spatial and temporal grid parameters $h$ and
$\tau$ defining $s_j^n \approx s(x_j,t_n)$, with $x_j=jh, j=0,\pm
1,\pm 2,\ldots$ and $t_n=n\tau,n =0,1,\ldots$. The vector
$S_h^n=\left(s_j^n\right)_{j=-\infty}^{j=\infty}$, and its grid 2-norm
$\norm{S_h^n}{}= \sqrt{h \sum_j (s_j^n)^2}$ is equivalent to and
indistinguished from the norm in $L^2(\R)$.  The discrete Fourier
transform applied to $S_h^n=(s_j^n)_j$ gives
$\widehat{S}^n=(\widehat{s^n}(\xi))_{\xi}$, with $-\pi/h \leq \xi\leq
\pi/h$. By Parseval's relation, the study of the evolution of
$\norm{S_h^n}$ is equivalent to the study of $\norm{\widehat{S^n}}$
defined through $L^2(-\frac{\pi}{h},\frac{\pi}{h})$. For one step
scheme from $t_{n-1} \to t_n$ we derive a formula for
\ba
\label{eq:gfactor}
\widehat{s}^n(\xi)=g(\xi)\widehat{s}^{n-1}(\xi).
\ea
The amplification factor $g(\xi)$ for $L=D$ and $L=A$ is well known;
see Tab.~\ref{tab:schemes} for the concise summary of conditions
required to establish a bound $\abs{g(\xi)}\leq 1$, from which the
strong stability $\norm{S^n_h} \leq \norm{S^{n-1}_h} $ follows, the
scheme is strongly stable, and the error propagates linearly. 

\begin{table}
\begin{center}
\begin{tabular}{l||c|c}
\hline\hline
Problem, scheme&$g(\xi)$&Condition\\
\hline
diffusion $L=D$, explicit $n*=n-1$&
$g_{DE}=1-s_D(\xi)$
&$4D\frac{\tau}{h^2}\leq 2$ 
\\
diffusion $L=D$, implicit $n*=n$&
$g_{DI} = (1+s_D(\xi))^{-1}$
&
none
\\
advection $L=A$, explicit $n*=n-1$&
$g_{AE}=1-s_{\lambda}$
&$q\geq 0, \tau\leq \frac{h}{q}$
\\
advection $L=A$, implicit $n*=n$&
$g_{AI}=(1+s_{\lambda})^{-1}$
&
$q \geq 0$
\\
\hline
\hline
\end{tabular}
\end{center}
\caption{\label{tab:schemes} Amplification factors and stability
  conditions for scalar diffusion and upwind advection schemes.
  Here $s_D(\xi)=2D\frac{\tau}{h^2}(1-\cos(\xi h))$,
  $\lambda=q\frac{\tau}{h}$, and $s_{\lambda}=\lambda(1-e^{-i\xi h})$.
}
\end{table}

\paragraph{\bf Von-Neumann analysis for systems} 
When approximating 
\ba
w(x,t)=[u(x,t),v(x,t)]^T \in \R^2,
\ea
the Fourier analysis is applied to each component of $w(x,t)$. For the
evolution system considered in this paper, instead of
\eqref{eq:gfactor}, we derive the system
\ba
\label{eq:h0h1}
	H_1 \widehat{w^{n}}(\xi)= H_0 \widehat{w^{n-1}}(\xi),
\ea
where $H_1,H_0 \in \C^{2 \times 2}$ are matrices dependent on $h,\tau$,
the Fourier variable $\xi$, and the coefficients of the
PDE.  The form directly resembling \eqref{eq:gfactor} is  
\ba
\label{eq:gmatrix}
	\widehat{w^{n}}= G(\xi) \widehat{w^{n-1}},
\ea
with the amplification matrix $G=G(h,\tau;\xi)={(H_1)}^{-1}H_0$. 

Our stability analysis establishes the conditions upon which the n'th
power $(G(\cdot;\xi))^n$ of $G$ is uniformly bounded for all $0\leq
n\tau \leq T$. For normal matrices, it suffices to study the spectral
radius $\rho(G)$, since, when $G^*G=GG^*$, all three members in the
inequality $ {\rho(G)}^n \leq \|G^n\| \leq {\|G\|}^n, n\geq 1, $
involving the 2-matrix norm ${\|G\|}$, are equal. For non-normal
matrices, one must study the spectral radius of $G^*G$ matrix, i.e.,
the first singular value of $G$, and the analysis of ${\|G\|}$ gets
quickly quite complicated. The symmetrizaton and change of variables
help in the calculations which otherwise are not easily manageable. 

\begin{proposition}
\label{prop:vonneumann}
Let $k=1$ and $q\geq 0, d \geq 0$. The fully discrete finite
difference schemes for \eqref{eq:bdiscrete} are strongly stable in
$[\sqrt{c}u,v]^T$ variables under the same conditions that apply to
the scalar diffusion, advection shown in Table~\ref{tab:schemes}. In
particular, the implicit schemes are unconditionally strongly stable
for $L=D$ and $L=A$, and the explicit-implicit schemes are
conditionally strongly stable for each $L=D$ and $L=A$.  The scheme
for $L=D+A$ in which diffusion is implicit, and advection is explicit
is strongly stable under the same conditions as that for explicit
advection.
\end{proposition}

The proof is established in the individual subsections below. 

\subsection{Stability of finite difference scheme for diffusion}
\label{sec:diffusion-ie}

We provide details for $L=D$ for the sake of exposition, since for a bounded
domain the case was already handled in Cor.~\ref{cor:iD} and
\ref{cor:eD} via MOL.

The row of a system \eqref{eq:bdiscrete-full} with $L_h$ as in
\eqref{eq:matrixD} is equivalent to 
\begin{subequations}
\label{eq:diff}
\begin{eqnarray}
	\frac{u_j^{n}-u_j^{n-1}}{\tau} + \frac{v_j^{n}- v_j^{n-1}}{\tau} + d\frac{2u_{j}^{n^*}-u_{j-1}^{n^*}-u_{j+1}^{n^*}}{h^2} & = & 0 \label{eq:diff-u}\\
	\frac{v_j^{n}- v_j^{n-1}}{\tau}+\alpha(v_j^{n}-cu_j^{n})&=&0 \label{eq:diff-v}
\end{eqnarray}
\end{subequations}
As suggested in Sec.~\ref{sec:abstract}, we first symmetrize
\eqref{eq:diff} by substituting \eqref{eq:diff-v} in
\eqref{eq:diff-u}, and rescale \eqref{eq:diff-u} by the factor
$\sqrt{c}$. Then we follow the usual non-Neumann analysis steps
applied to both components of $w_j^n=[\sqrt{c}u_j^n,v_j^n]^T$ and its
Fourier transforms
$[\sqrt{c}\widehat{u^n}(\xi),\widehat{v^n}(\xi)]^T$, which we denote
by $\widehat{w^n}$. (According to the convention adopted in
Sec.~\ref{sec:abstract} we should use $\widehat{\tilde{w}^n}$ but we will
  skip the tilde.)

\subsubsection{Implicit scheme for diffusion}
\label{sec:diffusion-i}

If $n*=n$, we rewrite \eqref{eq:diff} in the form \eqref{eq:h0h1} with 
\begin{equation}
	H_1=\left[ \begin{array}{c c} 1 +b c + s_D(\xi) & -b\sqrt{c}\\ -b\sqrt{c} & 1 + b \end{array}\right],\, H_0=I, \label{IMPdiffH1H0}
\end{equation} 
and where $b=\alpha \tau$. 

Since $s_D$ is real, thus $H_1$ is real and symmetric, and its
eigenevalues $\lambda_1,\lambda_2$ are real. Let $\lambda_1 \leq
\lambda_2$ while $\lambda_1+\lambda_2=Trace(H_1)=2+b(c+1)+s_D$. Since
$s_D\geq 0$, both the trace and the determinant of $H_1$ are positive,
with the latter given by
\bas
det(H_1)=1+b+bc+b^2c-b^2c+s_D(1+b).
\eas
We have thus
\bas
\lambda_1 \leq \frac{\lambda_1+\lambda_2}{2} \leq \lambda_2=\rho(H_1)=\norm{H_1}, 
\eas
and we see that  
$ {\| H_1\|} \geq 1 + \dfrac{b(c+1)}{2} + \dfrac{s_D(\xi)}{2} \geq 1 \mbox { for all } \xi$.
Since $s_D\geq 0$, we get
\ba
\norm{G}={\| H_1^{-1} \|}_2 \leq \frac{2}{2+b(c+1)}<1,
\ea 
which completes this case. 

\subsubsection{Explicit diffusion}
\label{sec:diffusion-e}
In this case $n*=n-1$, and we rewrite \eqref{eq:diff} in the form \eqref{eq:h0h1} with 
\ba
	H_1=\left[ \begin{array}{c c} 1 +bc & -b\sqrt{c}\\ -b\sqrt{c} & 1 + b \end{array}\right], \;\;
	H_0=\left[ \begin{array}{c c} 1 -s_D (\xi)& 0\\ 0 & 1 \end{array}\right] 
\ea
and both $H_1,H_0$ are real symmetric matrices. 

First we want to find a lower bound for
$\norm{H_1}=\rho(H_1)$. Denoting the eigenvalues of $H_1$ by
$\lambda_1,\lambda_2$, we calculate that $det(H_1)=1+ b(1+c)=
\lambda_1 \lambda_2$, and $Trace(H_1)=2+ b(1+c)=
\lambda_1+\lambda_2$. From this we see $\lambda_2(\lambda_1-1)=
\lambda_1 -1$, thus either $\lambda_1$ or $\lambda_2$ must equal
  1. Assuming, wlog, $\lambda_1=1$, we conclude, from
  $1+b(1+c)=\lambda_2$ that $\lambda_2>1$, and thus ${\| H_1^{-1}
    \|}_2 = 1.$

On the other hand, the eigenvalues for $H_0$ are on its diagonal, thus
the spectral radius is $\rho(H_0)=\max\{1,\abs{1-s_D}\}$. In order to
  guarantee $\norm{G}=\norm{H_1^{-1}H_0}\leq \norm{H_1^{-1}}\norm{H_0}
  \leq \norm{H_0}=\rho(H_0)\leq 1$, we must therefore have that $s_D\leq 2$ which 
requires
\begin{equation}
	\frac{d\tau}{h^2}<\frac{1}{2}.	 \label{eq:exp-diff-cond}
\end{equation}
In summary, the scheme is strongly stable if \eqref{eq:exp-diff-cond}
holds. This is the same result as that in Cor.~\ref{cor:eD} obtained by MOL. 

\subsection{Scheme for advection}
\label{sec:advection-discrete}
We begin by writing the upwind advection scheme for $L=A$
\begin{subequations}
\label{eq:advection}
\begin{eqnarray}
	\frac{u_j^{n}-u_j^{n-1}}{\tau} + \frac{v_j^{n}- v_j^{n-1}}{\tau} 
+ q\frac{u_{j}^{n*}-u_{j-1}^{n*}}{h} & = & 0 \label{Sc01eq01}\\
	\frac{v_j^{n}- v_j^{n-1}}{\tau}+\alpha(v_j^{n}-cu_j^{n})&=&0 \label{Sc01eq02}
\end{eqnarray}
\end{subequations}
We proceed as in Section~\ref{sec:diffusion-ie}, with symmetrization
and rescaling, to determine the matrices $H_1$ and $H_0$ in
\eqref{eq:h0h1}.
\subsubsection{Explicit advection}
\label{sec:advection-e}
We find that when $n*=n-1$
\begin{equation}
H_1= \left[ \begin{array}{c c} 1+bc& -b\sqrt{c}\\ -b\sqrt{c}  & 1 + b \end{array}\right],\, H_0=\left[ \begin{array}{c c} 1-s_\lambda& 0 \\ 0 & 1  \end{array}\right].
\end{equation} 
Our analysis here is similar to that in Section~\ref{sec:diffusion-e}
from which we have $\rho(H_1^{-1})\leq 1$. We find that the system is
strongly stable provided $\rho(H_0)=\max(1,\abs{1-s_{\lambda}})\leq
1$, which holds provided $0 \leq \lambda \leq 1$, and requires
$q\geq0$ and $\tau \leq \frac{h}{q}$, i.e., the usual CFL condition.
\subsubsection{Implicit advection}
\label{sec:advection-i}
Intuitively, we expect to find unconditional stability for
  $n*=n$. Proceeding as in Section~\ref{sec:diffusion-ie} we find
that
\begin{equation}
	H_1=\left[ \begin{array}{c c} 1 +bc+s_\lambda(\xi) & -b\sqrt{c}\\ -b\sqrt{c} & 1 + b \end{array}\right], \;\;
	H_0=I_{2}.
\end{equation}
where $s_\lambda(\xi)$ is given as in Tab.~\ref{tab:schemes}, and has a
positive real part $Re(s_{\lambda})$.

Now $H_1$ is complex symmetric, but not normal, and this requires
extra work as compared to the cases before.  To show $\norm{H_1}\geq
1$ which will demonstrate unconditional stability, we need to estimate
the spectral radius of $K=H_1H_1^{*}$. Since $\sqrt{|det(K)|}\leq
\rho(K)$, if we prove that $\det(K)\geq 1$, we are done.

We first calculate $K$, simplifying some notation in $H_1=
\left[ \begin{array}{c c} X+iY & -\beta\\ -\beta &
    \gamma \end{array}\right]$ where we substituted
$X=1+bc+Re(s_{\lambda})$, $Y=Im(s_{\lambda})$, and $\beta=b\sqrt{c}$,
and $\gamma=1+b$. We get
\bas
K=\left[ \begin{array}{c c} X^2+Y^2+\beta^2 & 
-\beta(X+iY)-\beta \gamma\\ -\beta(X-iY)-\beta\gamma
    & \beta^2 +\gamma^2 \end{array}\right],
\eas

After some lengthy calculations and simplifications, we find that
$\det(K)=Y^2\gamma^2 +(X\gamma-\beta^2)^2$ which has a lower bound of
$(X\gamma-\beta^2)^2$. We can estimate this term from below, reverting
to the original constants in $H_1$ and using
$Re(s_{\lambda})=\lambda(1-cos(\xi h)) \geq 0$, to see that
\bas
X\gamma-\beta^2=1+bc+b+(1+b)\lambda(1-cos(\xi h)) \geq 1+bc+b \geq 1,
\eas
and we're done. 
	
\subsection{IMEX scheme for explicit advection and implicit diffusion}
\label{sec:diffusion:advection}
The discrete system is as follows
\begin{subequations}
\label{eq:advectiondiffusion}
\ba
\nonumber
\frac{u_j^{n}-u_j^{n-1}}{\tau} 
-\alpha ( v_j^{n}- cu_j^{n}) 
+ q\frac{u_{j}^{n-1}-u_{j-1}^{n-1}}{h}&\phantom{+}&
\\
+ \frac{d}{h^2} (-u_{j+1}^{n} +2u_j^{n} -u_{j-1}^{n}) &=&  0 \label{2ISc02eq01}\\
\frac{v_j^{n}- v_j^{n-1}}{\tau}
+\alpha (v_j^{n}-cu_j^{n})&=& 0 \label{2ISc02eq02}
\ea
\end{subequations}
We quickly see that the matrix $H_1$ is the same as in
Section~\ref{sec:diffusion-i} and the matrix $H_0$ is the same as in
Section~\ref{sec:advection-e}. The analysis in these sections
therefore gives the strong stability of the scheme provided the CFL
condition holds.

\paragraph{\bf Summary} With the last case we conclude the proof of 
Proposition~\ref{prop:vonneumann}. 

\section{Numerical examples}
\label{sec:simulations}

In this Section we illustrate the theory we developed in
Sec.~\ref{sec:discrete} with three examples. First, we show the lack
of stability of the discrete system in the natural product norm; our
simple example motivates the use of weighted spaces and
symmetrization.
Second, we consider an advection example for which we test the
convergence of the numerical scheme, and compare the solutions for
different $\alpha$ to the equilibrium case. Third, we provide an
example and convergence rates for diffusion.

\subsection{Instability in Euclidean norm on $\R$}
\label{sec:0d}

Here we let $H=\R$, with $L=0.1$, and $f=0$, and we consider a fully
implicit time discretization \eqref{eq:bdiscrete} of
\eqref{eq:bsys}. The initial condition $w^{0}=[1,1]^T$ is given.

The discrete system is solved for the approximations $w^{n}=
[u^{n},v^{n}]^T$ with a fully implicit scheme
\ba
\label{eq:iter}
w^{n}=(I+\tau \Bmat)^{-1} w^{n-1}, \; 
\ea
We use $c=5$, $\alpha=0.1$, and $\tau=0.2$. 

\begin{figure}
\caption{\label{fig:stability} Illustration of the lack of strong
  stability discussed in Sec.~\ref{sec:0d}. Left: the phase plot
  $(u,v)$ shows that the norm $\norm{w}_{H \times H}$ does not necessarily
  decrease. Right: the plot of the weighted norm $\norm{w}_{W_c}(t)$
decreases while $\norm{w}_{H \times H}(t)$ does not.}
  \includegraphics[width=.5\textwidth]{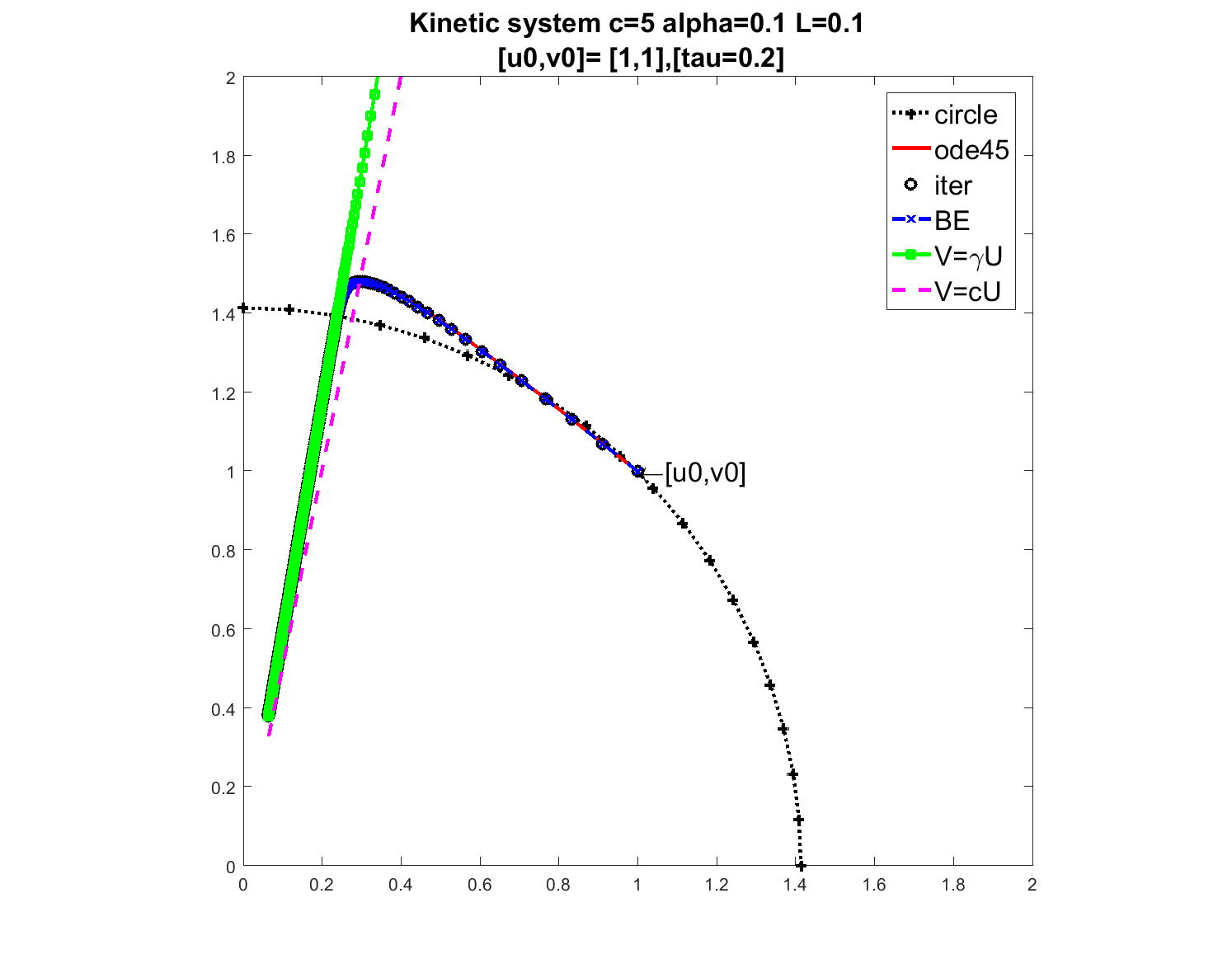}
  \includegraphics[width=.5\textwidth]{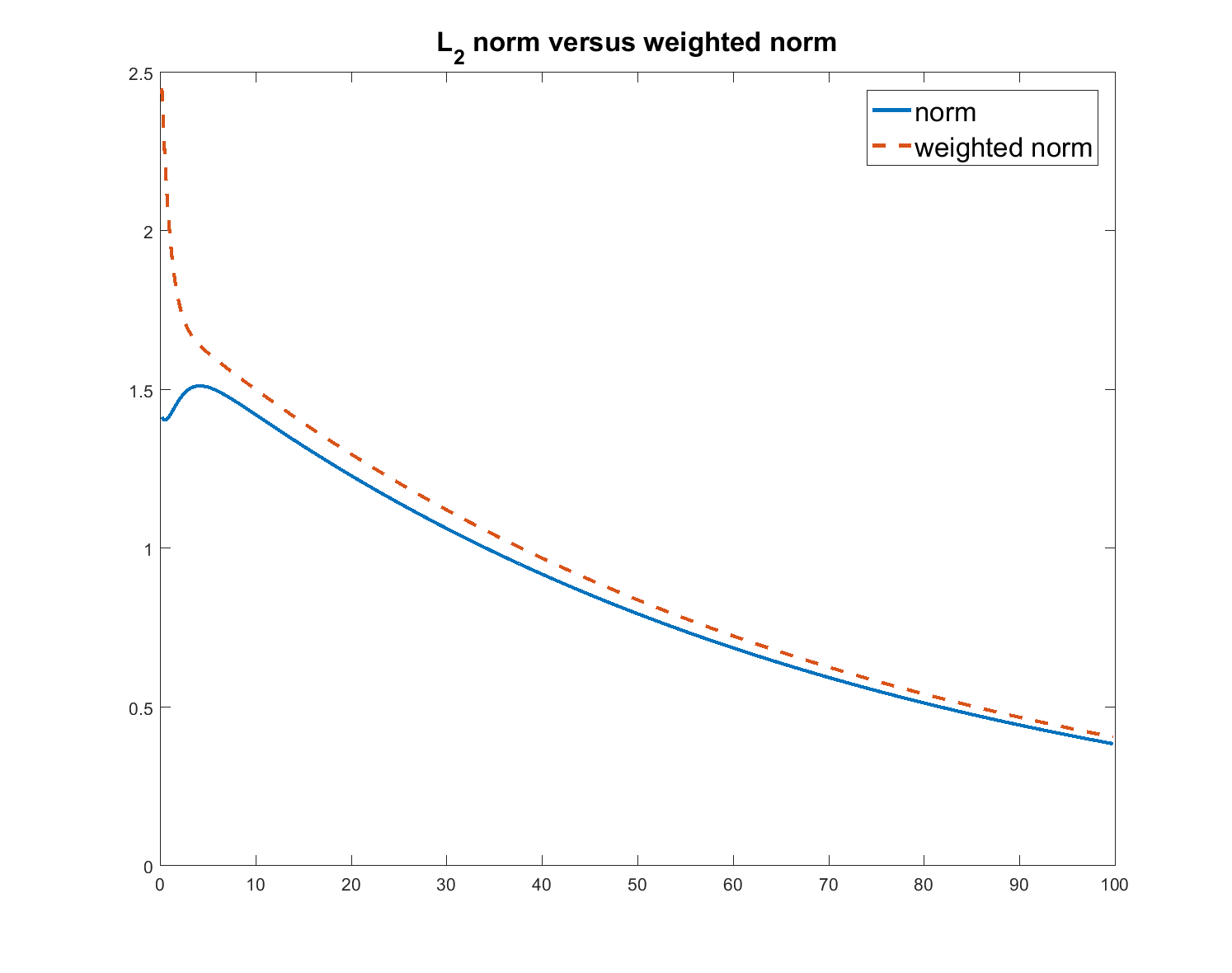}
\end{figure}

In Fig.~\ref{fig:stability} we illustrate the evolution of $w^n$;
these are close to those obtained to MATLAB's {\tt ode45} close to
$w^n$. It is clear that the solutions quickly tend to an asymptote and
then start decaying towards the origin. What is interesting is that,
the magnitude $w^n$ grows and the trajectory is ``above'' the circle
$\norm{w}=\norm{w^0}=\sqrt{2}$, before it heads towards the origin
along the asymptotic.

To explain, we examine $I+ \tau \Bmat$ which is not normal when $c \neq
1$.  In fact, even though its eigenvalues can be proven to be greater
than 1, its singular values are not both greater than 1. For example,
$\norm{(I+ \tau \Bmat)^{-1}}$ is $\approx 1.00741$, even though its
largest eigenvalue is $\approx 0.9971$. 

For independent interest, we study the asymptotics. To determine the
asymptotics, we solve for $v^{n}$ in terms of $u^{n}$, and substitute
back to \eqref{eq:bdiscrete}. Taking limits of both sides proves that
the limit, if it exists, must be $0$.  For the continuous problem
$w'+Bw=0$ we clearly expect that close to $[0,0]^T$ we will have $v$
follow close to $v=cu$. However, we find that $w^n$ actually follows
rather the asymptotics for the discrete system, $v \approx \gamma
u$. We can calculate the slope $\gamma$ directly from
\bas 
[u,\gamma u]^T=(I+\tau \Bmat)^{-1}[u,\gamma u]^T
\eas
In Figure~\ref{fig:stability} we illustrate both lines $v=cu$ and
$v=\gamma u$.

On the other hand, after symmetrization, the matrix $(I+\tau \tB)$ is
symmetric positive definite. We can calculate directly the eigenvalues
of $\tB=\tC+\Amat$, or simply show that for this symmetric $2 \times
2$ matrix, $\det(\tB)>0$ thus both of its eigenvalues
$\lambda_{1,2}\geq 0$. From this we conclude that the eigenvalues of
$I+\tau \tB$ are given by $1+\tau \lambda_{1,2}\geq 1$ and thus
$\norm{(I+\tau \tB)^{-1}}<1$. (For the numerical example as above, we
find $\norm{(I+\tau \tB)^{-1}} \approx 0.99732002$).

For illustration, we show that $\norm{w^{n}}_c$ is a decreasing
sequence but $\norm{w^n}$ is not.  This is illustrated also in
Figure~\ref{fig:stability}.

\subsection{Convergence of the schemes for advection and for diffusion}

With the stability results developed above, we expect the error for
the case $L=A$ to be of first order, as long as the true solution is
smooth enough. While the study of the regularity of the solutions is
outside the scope of this paper, we see that the case $L=A$ with
Riemann data develops enough smoothness to warrant first order error
in all $L_p$ spaces $1 \leq p < \infty$ and even for
$p=\infty$, similarly to what was observed in \cite{P13NumPDE}. In
turn, for $L=D$, with optimal smoothness, we expect second order
convergence, which is confirmed.

To test convergence, we use fine grid solution $u_{h_{fine}}$ instead of
manufacturing solutions which would require nonhomogeneous right-hand
side in \eqref{eq:ade2}. To simplify matters, we only report on
convergence rate at a fixed stopping time $T$. 

We define the error quantities (classical, and new quantity of interest)
\ba
E_{CQ}=\sqrt{{\| u-u_h \|}_{L_2}^2 +{\| v-v_h \|}_{L_2}^2},	
\\
E_{QoI}=\sqrt{c{\| u-u_h \|}_{L_2}^2 +{\| v-v_h \|}_{L_2}^2},	
\ea
where the $L_p$ grid norm for $1 \leq p < \infty$ is defined , as usual
\begin{equation}
	{\|u-u_h\|}_{L_p}= (\sum_i h {|u(x_i,T)-u_h(x_i,T)|}^p)^{1/p}.	
\end{equation}
In tables below, we report on the errors in different
  quantities of interest $E_r$ as well is in different norms
  $\norm{\cdot}{p}$, and calculate the respestive orders of the error
  $\alpha_r$, $\alpha_p$.

\subsubsection{Advection case}
We consider the problem 
\bas
u_t + v_t +u_x &=& 0, x \in \R\\
 v_t + \alpha(v-cu)&=&0.
\eas
and its approximation by the upwind scheme \eqref{eq:advection}. To
satisfy the CFL condition, we use $\lambda=0.99$, and we vary $\tau$
with $h$ in convergence terst. We choose intial data
\begin{subequations}
\begin{equation}
u(x,o)=	\mbox{``box''}(x)=\begin{cases} 1 & \mbox{ if } x \in [-1,0]\\ 0 &\mbox{ otherwise} \end{cases}.
\end{equation}
and $c=0.1,\,\alpha=2$. We also set 
\begin{equation}
v(x,0)=cu(x,0).
\end{equation}
\end{subequations}
which coresponds to \eqref{eq:initeq}. This helps to relate our
convergence rates to those obtained in \cite{P13NumPDE}. 

Since the true solution is not known, we use $M_{fine}=5050$ and
$T=4.8$. In Table~\ref{tab:advection} we show that the error
in every quantity of interest is of first order.

\begin{table}
\centering	
\begin{tabular}{cccccccc} 
M & ${\|u-u_h\|}_{L_2}$  & $\alpha_{L_2,u}$ & ${\|v-v_h\|}_{L_2}$ & $\alpha_{L_2,v}$ & ${\|u-u_h\|}_{L_1}$ & $\alpha_{L_1,u}$ \\ 
\hline	
20& 0.03682 & - & 0.004641 & - & 0.06217 & -\\
50& 0.01655 & 0.8728 & 0.002557 & 0.6504 & 0.02607 & 0.9483\\ 
100& 0.007575 & 1.127 & 0.0009912 & 1.367 & 0.01281 & 1.026\\ 
200& 0.003687 & 1.039 & 0.0004855 & 1.03 & 0.006244 & 1.036\\
500& 0.001329 & 1.113 & 0.0001771 & 1.101 & 0.002254 & 1.112\\
\hline
\end{tabular}
\centering	
\begin{tabular}{ccccccc} 	
M & ${\|u-u_h\|}_{\infty}$  & $\alpha_{\inf}$ & $E_{CQ}$ & $\alpha_{CQ}$ & $E_{QoI}$ & $\alpha_{QoI}$ \\ 
\hline		
20 & 0.03396 & - & 0.03711 & - & 0.1221 & -\\ 
50 & 0.02598 & 0.2922 & 0.01674 & 0.8685 & 0.05488 & 0.8728\\ 
100 & 0.007129 & 1.866 & 0.00764 & 1.132 & 0.02512 & 1.127\\ 
200 & 0.003529 & 1.014 & 0.003719 & 1.038 & 0.01223 & 1.039\\ 
500 & 0.001331 & 1.064 & 0.001341 & 1.113 & 0.004409 & 1.113 \\
\hline
\end{tabular}
\caption{Errors for advection case, with parameters $c=0.1,\,\alpha=2$, and ``box'' as the initial condition. Here $M_{fine}=5050$ and $T=4.8$}
\label{tab:advection}
\end{table}
 
\subsubsection{Convergence for diffusion}

We consider the problem 
\begin{subequations}
\bas
u_t + v_t -u_{xx} &=& 0, x \in (0,1)\\
 v_t + \alpha(v-cu)&=&0.
\eas
\end{subequations}
with the homogenous Dirichlet boundary conditions, and initial conditions
\begin{subequations}
\ba
u(x,0)=	\mbox{``bell''}= \exp\left(-{(x-0.5)}^2/0.3\right).
\label{inbell}		
\\
v(x,0)=cu(x,0).
\ea
\end{subequations}
In all experiments for this case we use $d=2$, stopping time $T=3.2$
and $M_{fine}=2000$. We vary $\tau=O(h^2)$ and expect optimal second
order convergence.
Indeed, Table~\ref{tab:diff} shows that error is $O(h^2)$ in
  every quantity of interest.
\begin{table}
\centering	
\begin{tabular}{ccccccc}M & ${\|u-u_h\|}_{L_2}$  & $\alpha_{L_2,u}$ & ${\|v-v_h\|}_{L_2}$ & $\alpha_{L_2,v}$ & ${\|u-u_h\|}_{L_1}$ & $\alpha_{L_1,u}$ \\
\hline	
20& 0.001785 & - & 0.01447 & - & 0.002758 & -\\ 
50 & 0.000285 & 2.003 & 0.002308 & 2.003 & 0.0004424 & 1.997\\
100 & 5.002$\cdot 10^{-5}$ & 2.51 & 0.0003468 & 2.734 & 8.73$\cdot 10^{-5}$ & 2.341\\
200 & 1.253$\cdot 10^{-5}$ & 1.997 & 8.657$\cdot 10^{-5}$ & 2.002 & 2.195$\cdot 10^{-5}$ & 1.992\\
\hline 
 \end{tabular}
\centering	
\begin{tabular}{ccccccc}M & ${\|u-u_h\|}_{\infty}$  & $\alpha_{\inf}$ & $E_{CQ}$ & $\alpha_{CQ}$ & $E_{QoI}$ & $\alpha_{QoI}$ \\ 
\hline	
 20.0 & 0.001849 & - & 0.01458 & - & 0.004373 & -\\
50.0 & 0.0002954 & 2.002 & 0.002326 & 2.003 & 0.000698 & 2.003\\
100.0 & 4.274$\cdot 10^{-5}$ & 2.789 & 0.0003504 & 2.73 & 0.0001225 & 2.51\\
200.0 & 1.067$\cdot 10^{-5}$ & 2.002 & 8.748$\cdot 10^{-5}$ & 2.002 & 3.07$\cdot 10^{-5}$ & 1.997\\
\hline
\end{tabular}
\caption{Errors for diffusion case, with parameters $c=5,\,\alpha=1.2$, and the ``bell'' as the initial condition. Here $M_{fine}=2000$; $T=3.2$}
\label{tab:diff}	
\end{table}

\subsection{Illustration of equilibrium vs non-equilibrium models}
\label{sec:noneq}
Now we are ready to show simulation results which illustrate the
kinetic effects in contrast to the equilibrium case. They are most
interesting for $L=A$, similar to \eqref{eq:advection}. We use
$\Omega=(-1,3)$ and periodic boundary conditions for
$u$.  We set $\alpha=2$ and $c=0.5$.

\begin{figure}
\caption{\label{fig:kinetic}Evolution for $L=A$ described in Sec.~\ref{sec:noneq}.}
\begin{center}
\includegraphics[width=.5\textwidth]{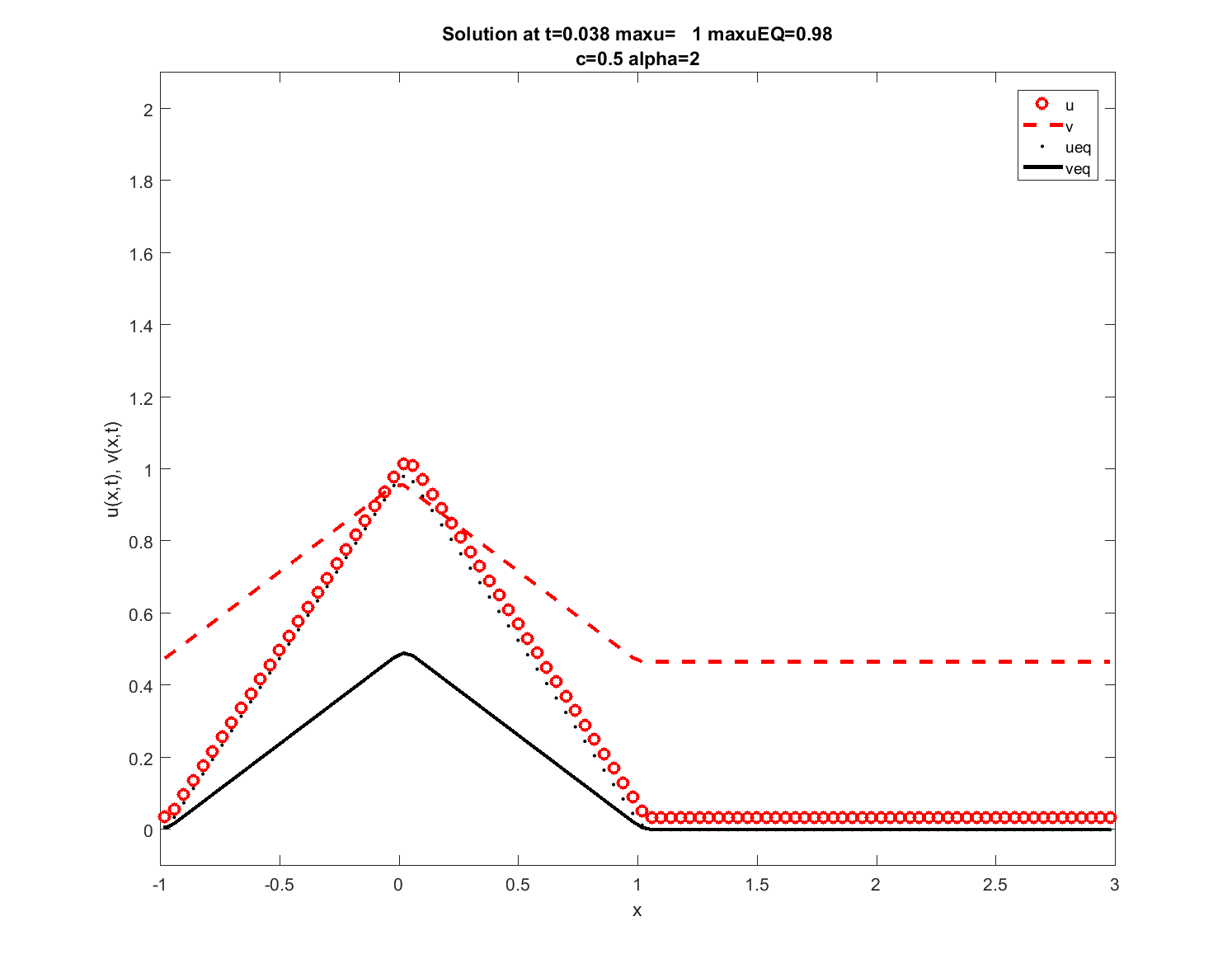}
\\
\includegraphics[width=.5\textwidth]{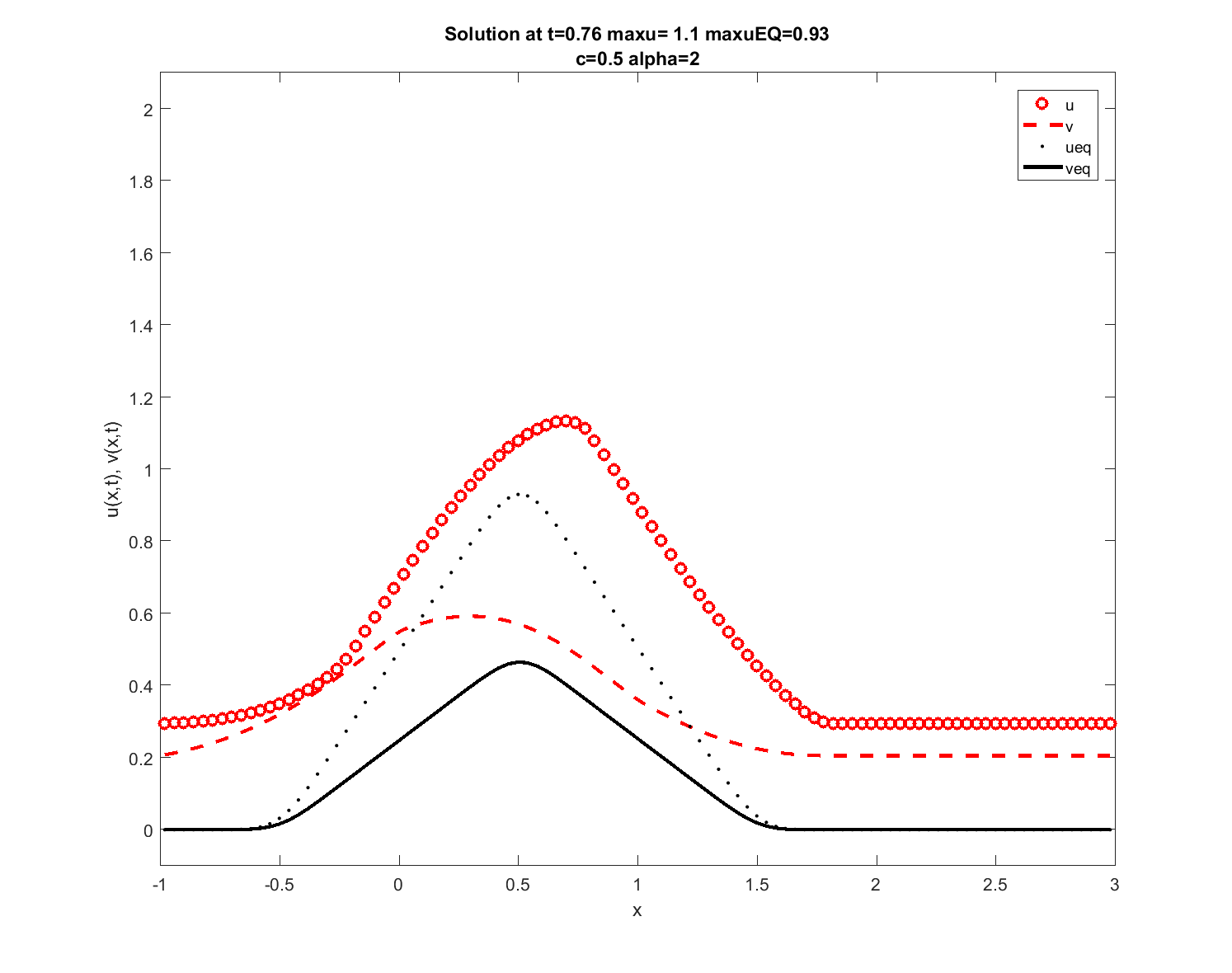}
\\
\includegraphics[width=.5\textwidth]{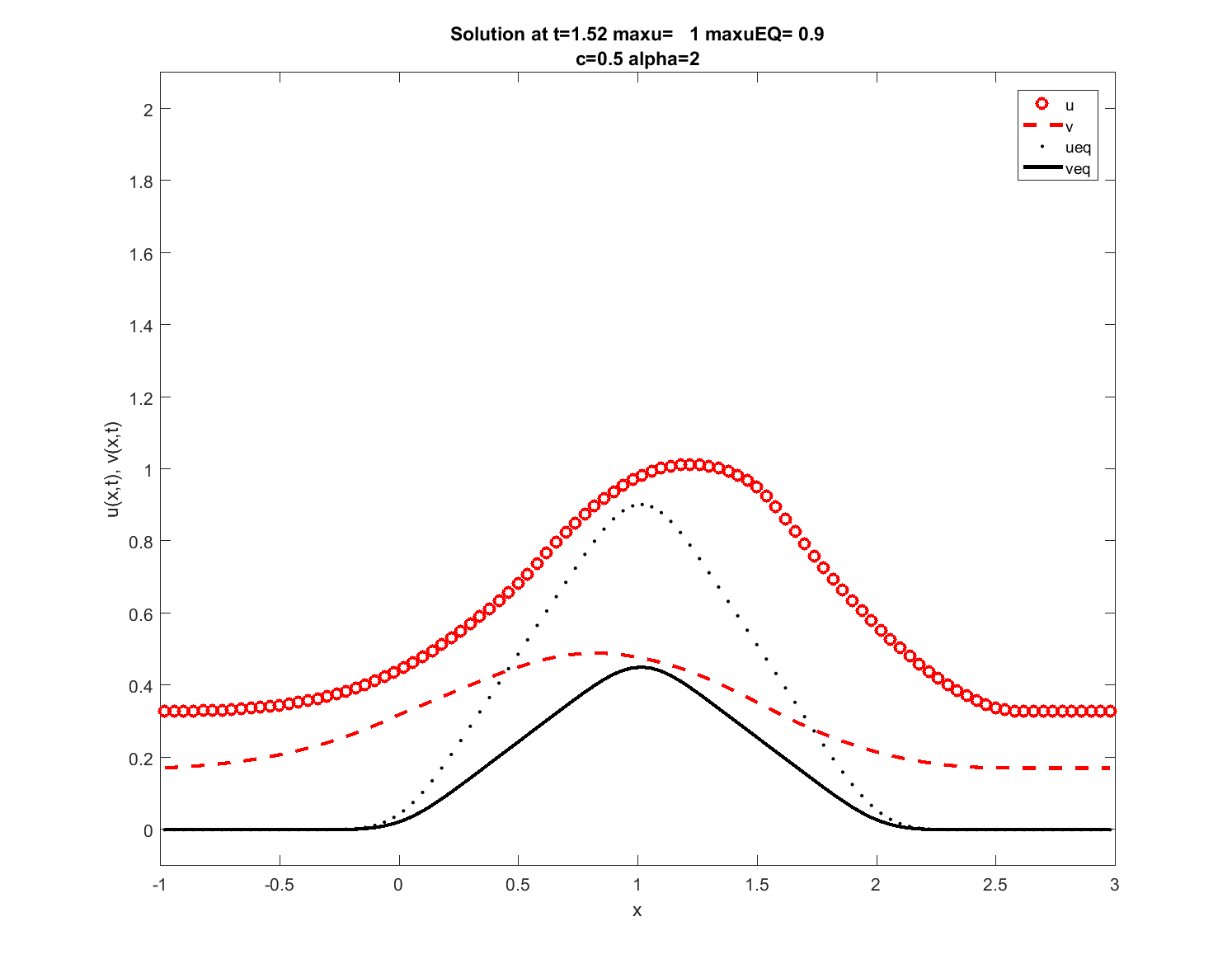}
\end{center}
\end{figure}

Fig.~\ref{fig:kinetic} shows the evolution of $[u(x,t),v(x,t)]^T$ at
  three time steps as shown. In addition, we show the evolution of an
  equilbrium solution in which $\alpha \to \infty$.

\section{Extensions}
\label{sec:extensions}

Above we have shown a unified framework for the analysis of
explicit-implicit schemes for the kinetic problems with a linear
non-equilibrium relationship. Further work is underway to generalize
these results, see below for nonlinear systems and systems with
multiple immobile sites. Error analysis and stability for
time-discrete schemes is underway. 

\subsection{Nonlinear equilibrium}
Here we consider the nonlinear extension of \eqref{eq:ade0} in which
$cU$ in \eqref{eq:ade02} is replaced by $g(u)$, with a monotone
increasing function $g:\R\mapsto R$. We only consider a finite
dimensional case and $H=\R^P$ since the proper setup with nonlinearity
in, e.g., $L^2(\Omega)$ is extensive and outside the present
scope. Here $g(U)=(g(u_j))_j, u \in \R^P$ is understood pointwise. The problem is 
\begin{subequations}
\label{eq:nde}
\ba
U'-\alpha(V-g(U))+LU&=& 0 \label{eq:nde1}
\\
V'+\alpha(V-g(U))&=& 0, \label{eq:nde2} 
\ea
\end{subequations}
and we will show stability of a particular new quantity of interest.

\begin{lemma}
\label{lem:gA}
Let $L$ and $g$ satisfy 
\ba
\label{eq:gA}
\langle LU,g(U) \rangle \geq 0
\ea
Then it holds that 
\ba
\frac{d}{dt}\left[ G(U) + \frac{1}{2}\norm{V}^2 \right] \leq 0.
\ea 
where $G$ is the primitive of $g(\cdot)$.
\end{lemma}
{\noindent \sc Proof: }To show the stability, we take the inner product of
\eqref{eq:nde1} with $g(U)$ and of \eqref{eq:nde2} with $V$
\begin{eqnarray*}
\langle U',g(U) \rangle -\alpha \langle V,g(U) \rangle
+\alpha \langle g(U),g(U) \rangle +
\langle AU,g(U) \rangle &=& 0 
\\
\langle V',V \rangle +\alpha \langle V,V \rangle -\alpha \langle g(U),V \rangle &=&
        0.
\end{eqnarray*}
Adding the two equations up we have 
\bas
\langle U',g(U)\rangle + 
\langle V',V \rangle 
+\alpha \langle V,V \rangle-2 \alpha \langle V,g(U) \rangle 
+\alpha \langle g(U),g(U) \rangle +
\langle AU,g(U)\rangle =0.
\eas
Rewriting we obtain 
\bas
\langle U',g(U) \rangle + \langle V',V \rangle +\alpha \left[ 
\langle V-g(U),V-g(U) \rangle \right]+ \langle LU,g(U) \rangle= 0.
\eas
Next step is to define a primitive $G:\R\to R$ of $g(\cdot)$ so that $G'(r)=g(r)$
and $\frac{d}{dt} G(U_j)=g(u_j)\frac{d}{dt}u_j$. 
We can thus write $ \langle U',g(U) \rangle =\frac{d}{dt}\sum_j G(u_j)$. 

Thus, if \eqref{eq:gA} holds, we obtain stability since
and we have proven Lemma~\ref{lem:gA}.

Next we provide sufficient conditions for \eqref{eq:gA} to hold. 
By a difference matrix \cite{Strang} we mean $D:\R^N\mapsto R^{N+1}$
such that $(Du)_j=u_j-u_{j-1}$ for $j=2,\ldots N-1$, and $(Du)_1=u_1,
(Du)_{N+1}=-u_N$. $D$ is therefore a discrete analogue of the
derivative (gradient). In turn, $D^T:\R^{N=1}\mapsto \R^N$ which
satisfies $(DU,V)_{R^{N+1}}=(U,D^TV)_{R^N}$, for any $U \in R^{N}$ and
$V \in R^{N+1}$, is the discrete analogue of the negative of the
divergence, which is dual to the gradient. (The analogues make sense
if one also assumes that $u_0=0$ and $u_{N+1}=0$, i.e., imposes
homogeneous Dirichlet boundary conditions on $U$.)

\begin{proposition}
\label{prop:gA}
Let $L=D^TKD$ where $D$ is a difference matrix, and $K$ is a diagonal
matrix with positive entries. Then \eqref{eq:gA} holds.
\end{proposition}
{\noindent \sc Proof: }It remains to prove that $L=D^TKD$ satisfies \eqref{eq:gA}.  We
consider first the case when $K=I$. The matrix $D^TD$ is the
well known ``discrete Laplacian'' $L_h$ in Section~\ref{sec:D}
which is symmetric positive definite, and it is easy to see that
$(D^TDu,u)_{R^N} = (Du,Du)_{R^{N+1}}=
u_1^2+\sum_{j=2}^{N}(u_j-u_{j-1})^2+U_N^2>0$ unless $U=0$.  Similarly
we obtain $(D^TDU,g(U))_{R^N} = (DU,Dg(U))_{R^{N+1}}=
u_1g(u_1)+\sum_{j=2}^N(u_j-u_{j-1})(g(u)-g(u_{j-1})+u_Ng(u_N)$ which
is nonnegative by virtue of $g(\cdot)$ being a monotone
increasing function.

In the more general case when $K \neq I$ we see that the argument
above holds for diagonal matrix $K$ with the the entries $k_1,\ldots
k_{N+1}$. Then we obtain $(D^TKDu,g(u))_{R^N} =
k_1u_1g(u_1)+\sum_{j=2}^N(k_j)(u_j-u_{j-1})(g(u)-g(u_{j-1})+k_{N+1}u_Ng(u_N)$. Since
each of these entries is nonnegative, we obtain the desired result. 

\subsection{Stability for a system with two species}
Here we consider again the finite dimensional space $H = \R^n$ and write 
\begin{eqnarray}
	U'+V_1'+V_2'+LU&=& 0 \label{eq-multu}\\ 
V_1'+\alpha_1(V_1-c_1U)&=&0
\\
V_2'+\alpha_2(V_2-c_2U)&=&0
\end{eqnarray}
Here we take the inner product of the first equation with $c_1c_2 U$,
the second by $c_2V_1$, the third by $c_1V_2$, (notice the
crossmultiplications) and add up to get
\begin{eqnarray*}
c_2 V_1^TV_1 +  c_2 V_2^TV_2 + c_1 c_2 U^TU + c_1 c_2 LU^2 
\\
+ \alpha_1 c_2 (V_1 -c_1U)^2 + \alpha_2 c_1 (V_2 -c_2U)^2 =0
\end{eqnarray*}
from which the stability follows for the following quantity
\ba
\frac{d}{dt} \left(c_1 c_2 \norm{U}^2 + c_2 \norm{V_1}^2 + c_1 \norm{V_2}^2\right) \leq 0. 
\ea
Further extensions to $m$ species are possible but will not be discussed. 

\section{Acknowledgements}

The authors would like to thank the anonymous reviewers whose
  remarks helped to improve the paper.

Research presented in this paper was partially supported by NSF grants
DMS-1115827 ``Hybrid modeling in porous media'', and DMS-1522734
``Phase transitions in porous media across multiple scales''; second
author served as a Principal Investigator on these projects. Majority
of research was done when first author F.~Patricia Medina was a PhD
student and later a faculty at Oregon State University.
 
\bibliography{mpesz,peszynska,BibKinStabil}
\end{document}

%% file: amsnumbers.tex
\newcommand{\0}{\mathbf{0}}

\newcommand{\N}{\mathbb{N}}
\newcommand{\C}{\mathbb{C}}

\def\ba{\begin{eqnarray}}
\def\ea{\end{eqnarray}}
\def\bas{\begin{eqnarray*}}
\def\eas{\end{eqnarray*}}

\newcommand{\bex}[1]{\begin{example}{\bf #1}}
\newcommand{\eex}{\end{example}}
\newcommand{\bv}{\begin{verbatim}}
\newcommand{\ev}{\end{verbatim}}
\def\bi{\begin{itemize}}
\def\ei{\end{itemize}}

%% file: MedinaPeszynska_2017_revised.bbl
\def\cprime{$'$} \def\cprime{$'$}
\begin{thebibliography}{10}
\expandafter\ifx\csname url\endcsname\relax
  \def\url#1{\texttt{#1}}\fi
\expandafter\ifx\csname urlprefix\endcsname\relax\def\urlprefix{URL }\fi
\expandafter\ifx\csname href\endcsname\relax
  \def\href#1#2{#2} \def\path#1{#1}\fi

\bibitem{BKneq97}
J.~W. Barrett, P.~Knabner,
  \href{http://dx.doi.org/10.1137/S0036142993249024}{Finite element
  approximation of the transport of reactive solutes in porous media. {I}.
  {E}rror estimates for nonequilibrium adsorption processes}, SIAM J. Numer.
  Anal. 34~(1) (1997) 201--227.
\newblock \href {http://dx.doi.org/10.1137/S0036142993249024}
  {\path{doi:10.1137/S0036142993249024}}.
\newline\urlprefix\url{http://dx.doi.org/10.1137/S0036142993249024}

\bibitem{BKneq95}
J.~W. Barrett, H.~Kappmeier, P.~Knabner, Lagrange-{G}alerkin approximation for
  advection-dominated contaminant transport with nonlinear equilibrium or
  non-equilibrium adsorption, in: Modeling and computation in environmental
  sciences ({S}tuttgart, 1995), Vol.~59 of Notes Numer. Fluid Mech., Vieweg,
  Braunschweig, 1997, pp. 36--48.

\bibitem{DvDW94}
C.~N. Dawson, C.~J. van Duijn, M.~F. Wheeler, Characteristic-{G}alerkin methods
  for contaminant transport with nonequilibrium adsorption kinetics, SIAM J.
  Numer. Anal. 31~(4) (1994) 982--999.

\bibitem{STW97}
H.~J. Schroll, A.~Tveito, R.~Winther, An l1--error bound for a semi-implicit
  difference scheme applied to a stiff system of conservation laws, SIAM
  journal on numerical analysis 34~(3) (1997) 1152--1166.

\bibitem{P13NumPDE}
M.~Peszynska,
  \href{http://www.math.oregonstate.edu/~mpesz/documents/publications/P13NMPDE.pdf}{Numerical
  scheme for a conservation law with memory}, Numerical Methods for PDEs 30
  (2014) 239--264.
\newblock \href {http://dx.doi.org/10.1002/num.21806\&ArticleID=1159335}
  {\path{doi:10.1002/num.21806\&ArticleID=1159335}}.
\newline\urlprefix\url{http://www.math.oregonstate.edu/~mpesz/documents/publications/P13NMPDE.pdf}

\bibitem{PSY15}
M.~Peszynska, R.~Showalter, S.-Y. Yi,
  \href{http://www.math.ualberta.ca/ijnam/Volume-12-2015/No-3-15/2015-03-04.pdf}{Flow
  and transport when scales are not separated: Numerical analysis and
  simulations of micro- and macro-models}, International Journal Numerical
  Analysis and Modeling 12 (2015) 476--515.
\newline\urlprefix\url{http://www.math.ualberta.ca/ijnam/Volume-12-2015/No-3-15/2015-03-04.pdf}

\bibitem{P95c}
M.~Peszy{\'n}ska, Finite element approximation of diffusion equations with
  convolution terms, Math. Comp. 65~(215) (1996) 1019--1037.

\bibitem{Show-monotone}
R.~E. Showalter, Monotone operators in {B}anach space and nonlinear partial
  differential equations, Vol.~49 of Mathematical Surveys and Monographs,
  American Mathematical Society, Providence, RI, 1997.

\bibitem{ShiMazumder08}
J.-Q. Shi, S.~Mazumder, K.-H. Wolf, S.~Durucan,
  \href{http://dx.doi.org/10.1007/s11242-008-9209-9}{Competitive methane
  desorption by supercritical {CO2}; injection in coal}, Transport in Porous
  Media 75 (2008) 35--54, 10.1007/s11242-008-9209-9.
\newline\urlprefix\url{http://dx.doi.org/10.1007/s11242-008-9209-9}

\bibitem{Kovscek2007}
K.~Jessen, W.~Lin, A.~R. Kovscek, Multicomponent sorption modeling in {ECBM}
  displacement calculations, SPE 110258.

\bibitem{JessenKovc08}
K.~Jessen, G.~Tang, A.~R. Kovscek, Laboratory and simulation investigation of
  enhanced coalbed methane recovery by gas injection, Transport in Porous Media
  73 (2008) 141--159.

\bibitem{PIMA11}
M.~Peszynska, Methane in subsurface: mathematical modeling and computational
  challenges, in: C.~Dawson, M.~Gerritsen (Eds.), IMA Volumes in Mathematics
  and its Applications 156, Computational Challenges in the Geosciences,
  Springer, 2013.

\bibitem{Showalter77}
R.~E. Showalter, Hilbert space methods for partial differential equations,
  Electronic Monographs in Differential Equations, San Marcos, TX, 1994,
  electronic reprint of the 1977 original.

\bibitem{KingErtekin86}
G.~R. King, T.~Ertekin, F.~C. Schwerer, Numerical simulation of the transient
  behavior of coal-seam degasification wells, SPE Formation Evaluation 2 (1986)
  165--183.

\bibitem{ShiDur03}
J.~Shi, S.~Durucan, A bidisperse pore diffusion model for methane displacement
  desorption in coal by {CO2} injection, Fuel 82 (2003) 1219--1229.

\bibitem{PS98}
M.~Peszy{\'n}ska, R.~E. Showalter, A transport model with adsorption
  hysteresis, Differential Integral Equations 11~(2) (1998) 327--340.

\bibitem{DiBSho82}
E.~DiBenedetto, R.~E. Showalter, A pseudoparabolic variational inequality and
  {S}tefan problem, Nonlinear Anal. 6~(3) (1982) 279--291.

\bibitem{WarRoo63}
J.~E. Warren, P.~J. Root, The behavior of naturally fractured reservoirs, Soc.
  Petro. Eng. Jour. 3 (1963) 245--255.

\bibitem{BZK60}
G.~I. Barenblatt, I.~P. Zheltov, I.~N. Kochina, Basic concepts in the theory of
  seepage of homogeneous liquids in fissured rocks (strata), J. Appl. Math.
  Mech. 24 (1960) 1286--1303.

\bibitem{ADH}
T.~Arbogast, J.~Douglas, Jr., U.~Hornung, Derivation of the double porosity
  model of single phase flow via homogenization theory, SIAM J. Math. Anal.
  21~(4) (1990) 823--836.

\bibitem{Show93a}
R.~E. Showalter, Diffusion in a fissured medium with micro-structure, in: Free
  boundary problems in fluid flow with applications (Montreal, PQ, 1990), Vol.
  282 of Pitman Res. Notes Math. Ser., Longman Sci. Tech., Harlow, 1993, pp.
  136--141.

\bibitem{HornSho90}
U.~Hornung, R.~E. Showalter, Diffusion models for fractured media, J. Math.
  Anal. Appl. 147~(1) (1990) 69--80.

\bibitem{KP12}
V.~Klein, M.~Peszynska, Adaptive double-diffusion model and comparison to a
  highly heterogenous micro-model, Journal of Applied Mathematics 2012 (2012)
  Article ID 938727, 26 pages.
\newblock \href {http://dx.doi.org/10.1155/2012/938727}
  {\path{doi:10.1155/2012/938727}}.

\bibitem{PS07}
M.~Peszy{\'n}ska, R.~E. Showalter, Multiscale elliptic-parabolic systems for
  flow and transport, Electron. J. Diff. Equations 2007 (2007) No. 147, 30 pp.
  (electronic).

\bibitem{TveWin}
A.~Tveito, R.~Winther, \href{http://dx.doi.org/10.1137/S0036141094263755}{On
  the rate of convergence to equilibrium for a system of conservation laws with
  a relaxation term}, SIAM J. Math. Anal. 28~(1) (1997) 136--161.
\newblock \href {http://dx.doi.org/10.1137/S0036141094263755}
  {\path{doi:10.1137/S0036141094263755}}.
\newline\urlprefix\url{http://dx.doi.org/10.1137/S0036141094263755}

\bibitem{BohmShow85b}
M.~B{\"o}hm, R.~E. Showalter, A nonlinear pseudoparabolic diffusion equation,
  SIAM J. Math. Anal. 16~(5) (1985) 980--999.

\bibitem{Thomee}
V.~Thom{\'e}e, Galerkin finite element methods for parabolic problems, 2nd
  Edition, Vol.~25 of Springer Series in Computational Mathematics,
  Springer-Verlag, Berlin, 2006.

\bibitem{Quarteroni}
L.~Fatone, P.~Gervasio, A.~Quarteroni, Multimodels for incompressible flows, J.
  Math. Fluid Mech. 2~(2) (2000) 126--150.

\bibitem{LeVeque}
R.~J. LeVeque, Finite difference methods for ordinary and partial differential
  equations, Society for Industrial and Applied Mathematics (SIAM),
  Philadelphia, PA, 2007, steady-state and time-dependent problems.

\bibitem{MortonRicht}
R.~D. Richtmyer, K.~W. Morton, Difference methods for initial-value problems,
  Second edition. Interscience Tracts in Pure and Applied Mathematics, No. 4,
  Interscience Publishers John Wiley \& Sons, Inc., New York-London-Sydney,
  1967.

\bibitem{Strang}
G.~Strang, Introduction to applied mathematics, Vol.~16, Wellesley-Cambridge
  Press Wellesley, MA, 1986.

\end{thebibliography}
